\documentclass[a4paper]{amsart}

\usepackage{verbatim,amsmath,amssymb,enumitem,amsthm}
\usepackage[displaymath,pagewise]{lineno} %for line numbers switch*,
\usepackage{mathtools,xcolor}

\newtheorem{thm}{Theorem} 

\newtheorem{prop}[thm]{Proposition}
\newtheorem{lem}[thm]{Lemma}
\newtheorem{cor}[thm]{Corollary}

\newtheorem{rem}[thm]{Remark}

\theoremstyle{remark}
\newtheorem{remark}[thm]{Remark}        % Numbered along with thm
\newtheorem{example}[thm]{Example}        % Numbered along with thm

\theoremstyle{plain} 
\newcommand{\thistheoremname}{}
\newtheorem*{genericthm*}{\thistheoremname}
\newenvironment{namedthm*}[1]
{\renewcommand{\thistheoremname}{#1}%
	\begin{genericthm*}}
	{\end{genericthm*}}

\numberwithin{equation}{section} \numberwithin{thm}{section}

\newtheorem{definition}[thm]{Definition}

\newcommand*\patchAmsMathEnvironmentForLineno[1]{%
	\expandafter\let\csname old#1\expandafter\endcsname\csname #1\endcsname
	\expandafter\let\csname oldend#1\expandafter\endcsname\csname end#1\endcsname
	\renewenvironment{#1}%
	{\linenomath\csname old#1\endcsname}%
	{\csname oldend#1\endcsname\endlinenomath}}%
\newcommand*\patchBothAmsMathEnvironmentsForLineno[1]{%
	\patchAmsMathEnvironmentForLineno{#1}%
	\patchAmsMathEnvironmentForLineno{#1*}}%
\AtBeginDocument{%
	\patchBothAmsMathEnvironmentsForLineno{equation}%
	\patchBothAmsMathEnvironmentsForLineno{align}%
	\patchBothAmsMathEnvironmentsForLineno{flalign}%
	\patchBothAmsMathEnvironmentsForLineno{alignat}%
	\patchBothAmsMathEnvironmentsForLineno{gather}%
	\patchBothAmsMathEnvironmentsForLineno{multline}%
}

\newcommand{\vf}{\varphi}
\newcommand{\cal}{\mathcal}

\newcommand{\dist}{{\rm dist}\,}

\newcommand{\R}{{\mathbb R}}

\newcommand{\INt}{{\rm int}\,}
\newcommand{\ep}{\varepsilon}

\newcommand{\diam}{{\rm diam}\,}

\newcommand{\Tan}{{\rm Tan}\,}

\newcommand{\spa}{\operatorname{span}}
\newcommand{\eps}{\varepsilon}

\newcommand{\sM}{\cal{M}}

\newcommand{\Q}{{\mathbb Q}}

\newcommand{\inter}{\mathrm{int}}

\def\en{\mathbb N}
\def\er{\mathbb R}

\newcommand{\graph}{\operatorname{graph}}
\newcommand{\card}{\operatorname{card}}

\newcommand{\spn}{\operatorname{span}}
\newcommand{\D}{\cal D}

\def\halfsq{\hbox{\kern1pt\vrule height 7pt\vrule width6pt height 0.4pt depth0pt\kern1pt}}
\def\ihalfsq{\hbox{\kern1pt \vrule width6pt height 0.4pt depth0pt
		\vrule height 7pt \kern1pt}}

\begin{document}

\title{A characterization of sets in $\R^2$ with DC distance function}
\author{Du\v san Pokorn\'y}\author{Lud\v ek Zaj\'i\v cek}
\address{Charles University\\
	Faculty of Mathematics and Physics\\
	Sokolovsk\'a 83\\
	186 75 Praha 8\\
	Czech Republic}

\thanks{The research was supported by GA\v CR~18-11058S}

\begin{abstract}  
	We give a complete characterization of closed sets  $F \subset \R^2$ whose distance function $d_F\coloneqq  \dist(\cdot,F)$ is DC
	(i.e., is the difference of two convex functions on  
	$\R^2$).  
	 Using this characterization, a number of properties of such sets is proved.
\end{abstract}  

\email{dpokorny@karlin.mff.cuni.cz}
\email{zajicek@karlin.mff.cuni.cz}

\date{\today}
\maketitle

\section{Introduction}

The present article is a continuation of the article \cite{PZ1} which studies closed
 sets $F \subset \R^d$,  whose distance function $d_F\coloneqq  \dist(\cdot,F)$ is DC
	(i.e., is the difference of two convex functions on $\R^d$). So we first briefly
	 recall the motivation for our study and mention some results of \cite{PZ1}.
	
	It is  well-known   (see, e.g., \cite[p. 976]{BB}) that, for a closed $F \subset \R^d$,
  the function $(d_F)^2$ is always DC but $d_F$ need not be DC. The main motivation
	 for the paper \cite{PZ1} was the  question whether
  $d_F$ a DC function if $F \subset \R^d$ is a graph of a $DC$ function $g: \R^{d-1} \to \R$.
This question  naturally arises
 in the theory of WDC sets (see  \cite[Question 2, p. 829]{FPR} and \cite[10.4.3]{Fu2}).
Let us note  that WDC sets  form a substantial generalization of Federer's sets
 with positive reach and still admit the definition of curvature measures (see \cite{PR} or \cite{Fu2}) and $F$ as in the above question is a natural example of a WDC set in $\R^d$. 
The main result of \cite{PZ1} gives the affirmative answer to the above question in the case
 $d=2$, but the case $d>2$ remains open. 

Following \cite{PZ1}, we will use the following notation.

\begin{definition}\label{D2}
For $d \in \en$, we set
$$ \cal D_d\coloneqq  \{\varnothing\} \cup  \{\varnothing \neq A \subset \R^d: \ A\ \ \text{is closed}\ \ \text{and}\ \  
d_A\ \ \text{is DC}\}. $$ 
The elements of  $\cal D_d$ will be called  $\cal D_d$ sets.            
\end{definition}

Using this notation, the main result of \cite{PZ1} asserts that
\begin{equation}\label{hlpz1}
\graph g \in \cal D_2,\ \text{whenever}\ g: \R \to \R \ \text{is DC}.
\end{equation}

If $A \subset \R^d$ is a set with positive reach, then (see \cite[Proposition 4.2]{PZ1})
 $A \in \cal D_d$ and also $\partial A \in \cal D_d$ and $\overline{\R^d \setminus A}\in \cal D_d$.
It implies (see \cite[Corollary 4.5]{PZ1}) that $\graph g \in \cal D_2$ whenever $g: \R^{d-1} \to \R$
 is a semiconcave function.

 It is not known whether each WDC set $A \subset \R^d$ belongs to $\cal D_d$, but the statement
 is true for $d=2$ by \cite[Theorem 3.3]{PZ2}. 
 
Several results concerning general properties of classes $\cal D_d$ were obtained in \cite[Section 4]{PZ1}; we recall
 them in subsection \ref{zname} below. 

In the present article, we use the results of \cite{PZ1} to give  complete characterizations
 of  $\cal D_2$  sets. These characterizations are based
 on the notion of (s)-sets (``special   $\D_2$-sets'') 
 which
 have a formally simple definition (Definition \ref{Ch})
 but their structure can be rather complicated. 
 The proofs of these characterizations are quite long and technical;
 they are contained in Sections 3 and 4. 

Section 5 contains applications of our characterizations. For example,
 we prove (Proposition \ref{sjdcgr}) that each nowhere dense $\D_2$ set is a countable
 union of DC graphs (defined in Definition \ref{dcgr}).
 Further, the system of all components of each $\D_2$ set is discrete
(see Theorem \ref{comp}) and each component is pathwise connected and locally connected (see
 Proposition \ref{souv}).
 An important application is  Theorem \ref{obrd2}; its particular case asserts that if $F: \R^2 \to \R^2$ is a bilipschitz
 bijection which is $C^2$ smooth (or, more generally, DC), then
 $F(M)\in \D_2$ for each $M\in \D_2$. It is an open question whether
 $\D_d$ has this stability property for $d>2$.

\section{Preliminaries}
\subsection{Basic notation}\label{basic}
We denote by $B(x,r)$ ($U(x,r)$) the closed (open) ball with centre $x$ and radius $r$.
 The boundary and the interior of a set $M$ are denoted by $\partial M$ and $\inter M$, respectively. A mapping is called $K$-Lipschitz if it is Lipschitz with a (not necessarily minimal) constant $K\geq 0$. In any vector space $V$, we use the symbol $0$ for the zero element and $\spa M$ for the linear span of a set $M$.

In the Euclidean space $\R^d$, the origin is denoted by $0$, the norm  by $|\cdot|$ and  the scalar product  by $\langle \cdot,\cdot\rangle$. By $S^{d-1}$ we denote the unit sphere in $\R^d$.  $\Tan(A,a)$ denotes the tangent cone of $A\subset \R^d$ at $a\in \R^d$ ($u \in \Tan(A,a)$ if and only if
 $u = \lim_{i \to \infty} r_i (a_i-a)$ for some $r_i>0$ and     $a_i \in A$,    $a_i \to a$).

If $x,y\in\R^d$, the symbol $\overline{x,y}$ denotes the closed segment (possibly degenerate). If also 
 $x \neq y$, then  $l(x,y)$  denotes the line
	 joining $x$ and $y$.
	 
	For $B\subset \er^2$ and $t \in \R$, we set  $B_{[t]}\coloneqq  \{y\in \R:\ (t,y) \in B\}$.    
We also define $\pi_1:\er^2\to\er$ by $\pi_1(x,y)=x$.

The distance function from a set  $\varnothing \neq A\subset \R^d$   is $d_A\coloneqq  \dist(\cdot,A)$ and 
 the metric projection of $z\in \R^d$ to $A$ is   
 $\Pi_A(z)\coloneqq \{ a\in A:\, \dist(z,A)=|z-a|\}$.

A system $\cal A$ of subsets of $\R^d$ is called discrete, if each
 $z\in \R^2$ has a neighbourhood which intersects at most one
 $A \in \cal A$. A set $D\subset \R^d$ is   called discrete,
  if $\{\{d\}: d\in D\}$ is discrete.

Under a rotation in $\R^2$ we always understand a rotation around the origin.

For  $f: \R^d \to \R^k$ and $x, v \in \R^d$, the one-sided directional derivative of $f$ at $x$ in direction $v$ is 
$$ f'_+(x,v)\coloneqq \lim_{t\to 0+} \frac{f(x+tv) - f(x)}{t}.$$

\subsection{DC functions}\label{dc}

 Let $f$ be a real function defined on an open convex set $C \subset \R^d$. Then we say
	 that $f$ is a {\it DC function}, if it is the difference of two convex functions.
	We say that $F=(F_1,\dots,  F_k): C \to \R^k$ is a {\it DC mapping} if all components
	 $F_i$ of $F$ are DC functions.

	 Semiconvex and semiconcave functions are special DC
	 functions. Namely, $f$ is a {\it semiconvex} (resp.
	 {\it semiconcave}) function, if  there exist $a>0$ and a convex function $g$ on $C$ such that
	$$   f(x)= g(x)- a |x|^2\ \ \ (\text{resp.}\ \  f(x)= a |x|^2 - g(x)),\ \ \ x \in C.$$

We will use the following well-known properties of DC functions and mappings.

\begin{lem}\label{vldc}
Let $C$ be an open convex subset of $\R^d$. Then the following assertions hold.
\begin{enumerate}
\item[(i)]
	If $f: C\to \R$ and $g: C\to \R$ are DC, then (for each $a\in \R$, $b\in \R$) the functions $|f|$, $af + bg$, 
	 $\max(f,g)$ and $\min(f,g)$ are DC.
	\item[(ii)]
	Each locally DC  mapping  $f:C \to \R^k$  is DC. 
		\item[(iii)]	Each  DC  function  $f:C \to \R$  is 
  Lipschitz on each compact convex set $Z\subset C$.
	%\item[(iv)] Let $f: C \to \R$ be a continuous function.
%Let $f_i: C \to \R$, $i=1,\dots,m$, be DC functions. Let $f: C \to \R$ be a continuous function
% such that $f(x) \in \{f_1(x),\dots,f_m(x)\}$ for each $x \in C$. Then $f$ is DC on $C$.
\item[(iv)] Let  $G\subset \R^d$ and $H\subset \R^k$ be open sets. Let $f: G \to \R^k$
  and $g:H \to \R^p$ be locally DC, and let $f(G) \subset H$. Then $g \circ f$ is locally DC
	 on $G$.
	\item[(v)]
	Let $G, H \subset \R^d$ be open sets, and let $f:G \to H$ be a  locally bilipschitz and locally
	 DC bijection. Then $f^{-1}$ is locally DC on $H$.
	\item[(vi)] Let $f_i: C \to \R$, $i=1,\dots,k$, be DC functions. Let $f: C \to \R$ be a continuous function
 such that $f(x) \in \{f_1(x),\dots,f_k(x)\}$ for each $x \in C$. Then $f$ is DC on $C$.
 	\item[(vii)] Each $C^2$ function $f:C\to\er$ is DC.
\end{enumerate}
\end{lem}
\begin{proof}
Property (i) follows easily from definitions, see e.g. \cite[p. 84]{Tuy}. Property (ii) was proved in \cite{H}.
 Property (iii) easily follows from the local Lipschitzness of  convex functions.
 Assertion (iv) is ``Hartman's superposition theorem'' from \cite{H}; for the proof see also
 \cite{Tuy} or \cite[Theorem 4.2]{VeZa}. Statement (v) follows from \cite[Theorem 5.2]{VeZa}.
 Assertion (vi) is a special case of \cite[Lemma 4.8.]{VeZa} (``Mixing lemma'').
   Property (vii) follows e.g. from \cite[Proposition 1.11]{VeZa} and (ii).  
\end{proof}
The following easy result (\cite[Lemma 2.3]{PRZ}) is well-known.
\begin{lem}\label{strict}
Let $F: (a,b) \to \R^d$ be a DC mapping and $x\in (a,b)$. Then the one-sided derivatives $F'_{\pm}(x)$ exist. Moreover
\begin{equation}\label{limder}
 \lim_{t\to x+} F'_{\pm}(t)=F'_+(x)\ \ \text{and}\ \       \lim_{t\to x-} F'_{\pm}(t)=F'_-(x)
\end{equation}
which implies that $F'_+(a)$ is the strict right derivative of $F$ at $x$, i.e.,
\begin{equation}\label{strd}
 \lim_{\substack{(y,z) \to (x,x)\\ y \neq z,\, y\geq x,\, z\geq x}}  \ \frac{F(z) - F(y)}{z-y} \ = \ F'_+(x).
 \end{equation}
\end{lem}

The notion of DC mappings between Euclidean spaces was generalized in \cite{VeZa} to the notion of
 DC mappings between Banach spaces  using the notion of a ``control function''.
%In the theory of DC mappings between Banach spaces (see \cite{VeZa}), it is essential the notion of %a control function for a DC mapping.
We will use this notion only for real functions defined on open intervals  $I \subset \R$. 
In this context we have 
 (cf. \cite[ Definition 1.1]{VeZa}) that a convex function $\vf:I \to \R$ is a {\it  control function} for a function $f: I \to \R$ if and only if both $\vf + f$ and $\vf -f$ are convex functions. 
 
 It is an easy fact (cf. \cite[Lemma 1.6 (b)]{VeZa})  that $f: I \to \R$ is DC if and only if it has a control function.
 We will use the following immediate consequence of  \cite[Proposition 1.13]{VeZa}.
 \begin{lem}\label{smercontr}
 If $\vf$ is a control function for $f$ on an open interval $I$, then
 $$  \left| \frac{f(z+k)-f(z)}{k} - \frac{f(z)-f(z-h)}{h}\right| \leq   \frac{\vf(z+k)-\vf(z)}{k} - \frac{\vf(z)-\vf(z-h)}{h},$$
  whenever $k>0$, $h>0$, $z\in I$, $z+k \in I$ and $z-h \in I$.
 \end{lem}
For the origin of the following definition, see \cite[p. 28]{RoVa}.
\begin{definition}
Let  $f$ be a function on $[a,b]$.
 For every partition
$D= \{a=x_0<x_1<\dots <x_n=b\}$ of $[a,b]$,  we put
$$
K (f,D)\coloneqq
\sum_{i=1}^{n-1} \left| \frac{f(x_{i+1})-f(x_i)}{x_{i+1} - x_i} -
\frac{f(x_{i})-f(x_{i-1})}{x_{i} - x_{i-1}} \right|.
$$
(If $n =1$, we put $K(f,D)\coloneqq  0$.) Then the {\em convexity of $f$ on $[a,b]$\/}
is defined as
$$ K_a^b f \coloneqq \sup  K (f,D),$$
where the supremum is taken over all partitions $D$ of $[a,b]$.
 If $K_a^b f< \infty$,
we say that {\em $f$ has a bounded (or finite) convexity}.
\end{definition}
The following fact is a consequence of \cite[Theorem 3.1(b)]{VZ2}.
\begin{lem}\label{ok}
%\begin{enumerate}
%\item[(i)]
%If $f$ is a continuous function on $[a,b]$, then $K_a^b f = \sup\{ %K_c^d f:\ a<c<d<b\}$.
%\item[(ii)]
If $f$ is a DC function on $(a,b)$ with a control function $\varphi$ and $a<c<d<b$, then
$K_c^d f \leq \varphi'_-(d)- \varphi'_+(c)$.
%\end{enumerate}
\end{lem} 
Following \cite[p. 1617]{PRZ}, we use the following terminology. 
\begin{definition}
We will say that a function defined on a set $\varnothing \neq D \subset \R$ is a {\it DCR  function} if it is a restriction of
 a DC function defined on $\R$.
\end{definition}

The following facts are well-known.
\begin{lem}\label{dcr}
 Let $f$ be a continuous real function on $[a,b]$. Then the following conditions are equivalent.
\begin{enumerate}
\item[(i)]
 $f$ is a DCR function.
\item[(ii)]
$f$ is the difference of two Lipschitz convex functions.
\item[(iii)]
$f$ has a bounded convexity.
\item[(iv)]
$f'_-(x)$ exists for each $x \in (a,b)$ and $V(f'_-, (a,b))< \infty$.
\item[(v)]
$f$ is a restriction of a DC function defined on some $(u,v) \supset
 [a,b]$.
\end{enumerate}
\end{lem}
(Here $V(f'_-, (a,b))$ means the variation of $f'_-$ over $(a,b)$ in the 
 usual sense; see, e.g. \cite[p. 322]{VZ2}.)
\begin{proof}
The implication $(i) \Rightarrow  (ii)$ follows by
 Lemma  \ref{vldc} (iii) and  
 $(ii) \Rightarrow  (i)$  holds since  each convex Lipschitz function on $[a,b]$
 can be  extended to a convex function on $\R$. The equivalence 
$(ii)\Leftrightarrow (iii)$ easily follows from \cite[Theorem~D, p. 26]{RoVa}
 and $(iii)\Leftrightarrow (iv)$ is a particular
 case of \cite[Proposition 3.4, p. 382]{VZ2}. The implication
  $(i) \Rightarrow (v)$  is trivial and $(v) \Rightarrow (iii)$
	 follows from Lemma \ref{ok}.
 \end{proof}

We will need the following facts concerning DCR functions. They immediately follow from \cite[Proposition 4.2]{VZ2} (or can be rather easily obtained using Lemma \ref{vldc} (iv),(v)).

\begin{lem}\label{vldcr}
Let $\vf:[a,b] \to [c,d]$ be a DCR increasing bilipschitz bijection
 and let $\omega: [c,d]\to \R$ be a DCR function. Then
\begin{enumerate}
\item[(i)] 
the function $\omega \circ \vf$ is DCR on $[a,b]$ and
\item[(ii)]
the function $\vf^{-1}$ is DCR on $[c,d]$.
\end{enumerate}
\end{lem}
We will need also the following    ``DCR mixing  lemma''.
 
\begin{lem}\label{L:DCRmixing}
	Let $I\subset \er$ be a closed interval and let $f: I \to \R$ be a continuous function.
	Let $f_i:I \to \R$, $i=1,\dots,k$, be DCR functions such that $f(x) \in \{f_1(x),\dots,f_k(x)\}$ for each $x \in I$. Then $f$ is DCR.
\end{lem}

\begin{proof}
	Let $\tilde f:\er\to\er$ be a continuous extension of $f$ which is locally constant on $\er\setminus I$ and let $\tilde f_i:\er\to\er$, $i=1,\dots,k$,
	 be a DC extension of $f_i$.
	 Then there are two constant (and so DC) functions $\tilde f_{k+1}$,
	 $\tilde f_{k+2}$ on $\R$ such that 
	  $\tilde f(x)\in\{\tilde f_1(x),\dots,\tilde f_{k+2}(x)\}$, $x \in \R$. Consequently $\tilde f$ is DC by Lemma~\ref{vldc}~(vi), and so
	 $f$ is DCR.	
	%	is a set $P$ of cardinality at most two such that $\tilde %f(x)\in\{\tilde f_1(x),\dots,\tilde f_k(x)\}\cup P$ and so 
	 %$\tilde f$ is DC by Lemma~\ref{vldc}~(vi) (note that constant %functions are DC).
	 %Since $f=\tilde f|_{I}$ we obtain that $f$ is DCR.
\end{proof}

The following lemma is a  version of the ``mixing lemma'' \cite[Lemma 4.8]{VeZa} (cf. 
 Lemma \ref{vldc} (vi)), which we need.
 Note that \cite[Lemma 4.8]{VeZa}  works even  with  DC mappings between Banach spaces, and 
  Lemma \ref{vmix} follows from its proof but not from
its
 formulation.

\begin{lem}\label{vmix}
Let $F_i$, $i=1,\dots,k$, be DC functions 
 on an open interval $J \subset \R$. Then there exists a convex function $\vf$ on $J$ with the following property:
\smallskip

(P)\ \ If $F$ is a continuous function on an open interval $I\subset J$ and $F(x) \in \{F_1(x),\dots,F_k(x)\}$, $x\in I$,
 then $F$ is a DC function with  control function $\vf|_I$.
\end{lem}
\begin{proof}
Let $f_i$ be a control function for $F_i$ on $J$, $i=1,\dots,k$. Set
$$ \vf\coloneqq  \sum_{i,j=1}^k\ h_{i,j},\ \ \ \text{where}\ \ \ h_{i,j}\coloneqq  
f_i + f_j + \frac 12 |F_i-F_j|.$$  
The proof of \cite[Lemma 4.8]{VeZa} (where $\vf$ is denoted
 by  $f$) gives  the assertion of property (P) for $I=J$. Observing
 that $f_i|_I$ is  a control function for $F_i|_I$, property
 (P) follows.
\end{proof}

We will need also the following easy ``Lipschitz mixing lemma''.

\begin{lem}\label{lipmix}
	Let  $K>0$ and   $f_i,\ i=1,\dots,k$, be $K$-Lipschitz functions on an interval (of arbitrary type)
	$I \subset \R$. Let $f$ be a continuous function on $I$ such that
	$f(x) \in \{f_1(x),\dots,f_k(x)\}$, $x \in I$. Then $f$ is $K$-Lipschitz on $I$.
\end{lem}
\begin{proof}
We will proceed by induction on $k$. The case $k=1$ is trivial. Suppose that $k>1$ and the lemma
 holds for ``$k=k-1$''. To prove that $f$ is $K$-Lipschitz, consider arbitrary points $a, b \in I$,
 $a<b$. Choose $1\leq i_0 \leq k$ such that $f(a)= f_{i_0}(a)$ and set
 $c\coloneqq  \max \{a\leq x \leq b:\ f(x)= f_{i_0}(x)\}$. 

If $c=b$, then $|f(b)-f(a)| = | f_{i_0}(b) - f_{i_0}(a)| \leq K (b-a)$.

If $c<b$, the induction hypothesis (applied to $f|_{(c,b]}$) implies that $f$ is
$K$-Lipschitz on $(c,b]$, and consequently also on $[c,b]$. Therefore
$$ |f(b)-f(a)| \leq | f_{i_0}(c) - f_{i_0}(a)| + |f(b)-f(c)| \leq K(c-a) + K (b-c) = K (b-a).$$
\end{proof}

By well-known properties of convex and concave functions, we easily obtain that each locally DC function $f$ on an open set $U \subset \R^d$ has all one-sided directional derivatives finite and
\begin{equation}\label{zlose}  
 g_+'(x,v) + g_+'(x,-v) \leq 0,\ \ x \in U, v\in \R^d,\ \ \ \text{if}\ \ g\ \ \text{is locally semiconcave on}\ \ U.
\end{equation}

%First note that the distance function $d_A$ and the  metric projection $\Pi_A$ were defined in
%Subsection~\ref{basic}.

Recall that if $\varnothing \neq A\subset \R^d$ is closed, then $d_A$ need not be DC; however (see, e.g., \cite[Proposition 2.2.2]{CS}), 
\begin{equation}\label{loksem}
\text{$d_A $ is locally semiconcave (and so locally DC) on $\R^d \setminus A$.}
\end{equation}

In \cite{PRZ} and \cite{PZ1} we worked with ``DC hypersurfaces'' in $\R^d$. Since we work here
 in $\R^2$ only, we use the following terminology.

\begin{definition}\label{surf}
We say that a set $A\subset \er^2$ is a $1$-dimensional DC surface, if there exist $v\in S^{1}$ and a Lipschitz DC function (i.e. the difference of two convex functions) $g$ on $W\coloneqq (\spn \{v\})^\bot$ such that $A=\{w+g(w)v: w\in  W\}$.
\end{definition}

\begin{rem}\label{hyper}
The notion of a $1$-dimensional DC surface in $\R^2$ 
coincides with the notion of a DC hypersurface in $\R^2$ from \cite{PRZ} (but not 
 with the notion of a DC hypersurface in $\R^2$ from \cite{PZ1}, where the Lipschitzness of
 $g$ is not required).
\end{rem}

We also define, following \cite{PZ1}, the  notion of a {\it DC graph} in $\R^2$.

\begin{definition}\label{dcgr}
A set $P\subset\er^2$ will be called a DC graph if it is a rotated copy of $\graph f$ of a DCR function $f$ on some compact (possibly degenerated) interval $\varnothing \neq I\subset\er$.
\end{definition}

Note that $P$ is a DC graph if and only if it is a nonempty connected compact subset of a $1$-dimensional DC surface in $\er^2$.

We will need the following simple result which is possibly new and can be of some independent interest.

\begin{prop}\label{obrdc}
	Let $g$ be a continuous function on $[a,b]$ which is DC on $(a,b)$ and let $P\subset [a,b]$ be a nowhere dense set.
	 Then the set $g(P)$ is nowhere dense.
\end{prop}

\begin{proof}
We can suppose that $P$ is closed. Suppose, to the contrary, that (the compact set) $g(P)$
 is not nowhere dense and choose an open interval $I \subset g(P)$. Set
$$ S\coloneqq  \{ x \in (a,b):\ g'_+(x)=0\ \ \text{or}\ \ g'_-(x)=0\}.$$
Then $g(S)$ is Lebesgue null; it follows e.g. from \cite[Theorem 4.5, p. 271]{Sa} (cf.
\cite[p. 272, a note before Theorem 4.7]{Sa}).
So we can choose a point $y_0 \in I \setminus(g(S)\cup\{g(a),g(b)\}).$
Then the set  $K\coloneqq  g^{-1}(\{y_0\}) \subset (a,b)$ is finite. Indeed, otherwise
 there exists a point $x\in K$ which is an accumulation point of the compact set $K$.
 Then clearly $x\in S$ which contadicts $y_0 \notin g(S)$. Let $K= \{x_1<...<x_p\}$.
   Lemma \ref{strict} implies that there exists $\delta>0$ such that
	$a<x_1-\delta<x_1+\delta< x_2-\delta<\dots< x_p+ \delta<b$ and $g$ is strictly monotone
	 both on $[x_i-\delta,x_i]$ and  on   $[x_i,x_i+\delta]$, $i=1,\dots,p$. Consequently
 	$Q\coloneqq g(P \cap \bigcup_{i=1}^p  (x_i-\delta, x_i + \delta))$  is nowhere dense.
	Since  $Z \coloneqq g( [a,b] \setminus  \bigcup_{i=1}^p  (x_i-\delta, x_i + \delta))$  is compact and does not
	 contain $y_0$, there exists $\sigma>0$ such that  $(y_0-\sigma, y_0+ \sigma) \subset I\subset g(P)$
	 and $ (y_0-\sigma, y_0+ \sigma)  \cap Z = \varnothing$.  Consequently $(y_0-\sigma, y_0+ \sigma)$ 
	  is a subset of nowhere dense set $Q$, which is a contradiction.
\end{proof}

\subsection{Known results concerning $\cal D_d$}\label{zname}

In, \cite{PZ1}, we proved several general results concerning systems $\cal D_d$. First recall
 that   if $M \subset \R$ is closed then
\begin{equation}\label{cald1}
\text{$M$ belongs to $\cal D_1$ iff the system of all components of $M$
 is locally finite.}
\end{equation}    
It easily implies that $\cal D_1$ is  closed with respect to both finite unions and finite
 intersections and that a closed $M\subset \R$ belongs to  $\cal D_1$ if and only if $\partial M \in \cal D_1$.

However, the case $d>1$ is different. It is easy to show that
\begin{equation}\label{sjedno}
\cal D_d\ \ \text{is closed with respect to finite unions}
\end{equation}
 but
\cite[Example 4.1]{PZ1} shows that already $\cal D_2$ is not closed with respect to finite intersections.
We observed that, for a closed set $M\subset \er^d$, $d\in\en$,
\begin{equation}\label{eq:boundaryEquivalenceGeneral}
\partial M \in \cal D_d\iff (M \in \cal D_d \quad\text{and}\quad \overline{\R^d \setminus M}\in  \cal D_d) 
\end{equation}
 but \cite[Example 4.1]{PZ1} provides an example of a set $M \in \D_2$ with $\partial M \notin \D_2$.

Important \cite[Proposition 4.7]{PZ1} asserts that if
 $d\geq 2$ and $M \in \cal D_d$, then each bounded set $C \subset \partial M$ can be covered by finitely many ``DC hypersurfaces''.

In our terminology we have in $\R^2$ the following result which has basic importance for the present article.
\begin{prop}\label{pokr}
Let $M \in \cal D_2$. Then each bounded set $C \subset \partial M$ can be covered by finitely many $1$-dimensional DC surfaces.
\end{prop}
Let us note that this result is not a particular case of \cite[Proposition 4.7]{PZ1} (cf. Remark
 \ref{hyper}), but the proof of \cite[Proposition 4.7]{PZ1} is based on \cite[Corollary 5.4]{PRZ}
 and so \cite[Proposition 4.7]{PZ1} holds also with the definition of DC hypersurfaces
 from \cite{PRZ} (``with Lipschitzness'') and so Proposition \ref{pokr} holds.
 (Moreover, it is easy to show that ``Proposition \ref{pokr} without Lipschitzness'' implies
 ``Proposition \ref{pokr} with Lipschitzness''.)

Proposition \ref{pokr} easily implies that each nowhere dense $M\in \cal D_2$ can be covered by 
 a locally finite system of DC graphs. On the other hand,  \eqref{hlpz1} easily implies (see \cite[Proposition 4.9]{PZ1})
  that 
	\begin{multline}\label{loksjgr}
	\text{if $M\subset\er^2$ is the union of a locally finite system}\\
	\text{of DC graphs then $M\in \cal D_2$.}
	\end{multline}

However, we have found an example (see \cite[Example 4.10]{PZ1}) of a nowhere dense set $M\in \cal D_2$ which
 is not the union of a locally finite system of DC graphs.

(Let us note that we will prove in Proposition \ref{sjdcgr} that each nowhere dense $M\in \cal D_2$ is the union
 of a countable system of DC graphs.)
 
 We will also use the following easy facts which are not mentioned explicitely in \cite{PZ1}.
 
\begin{rem}\label{rem:satbilityOfD2}
 	\begin{enumerate}
 		\item\label{cond:D2underSimilarity} If $M$ is a $\D_2$ set and $\varphi$ a similarity on $\er^2$ then $\varphi(M)$ is a $\D_2$ set.
 		\item\label{cond:D2locallyFiniteStable} The system $\D_2$ is closed with respect to locally finite unions.
 	\end{enumerate}
 \end{rem}
\begin{proof}
	The first part follows by   Lemma~\ref{vldc}~(i),~(ii),~(iv) and (vii)   from the fact that $d_{\varphi(M)}=rd_M\circ\varphi^{-1}$, where $r>0$ is the scaling ratio of $\varphi$.
	
	To prove the second part, consider a locally finite system $\sM\subset\D_2$ and denote $M\coloneqq\bigcup \sM$.
	   We can suppose that $M \neq \varnothing$. If $z \in M$, 
	choose $r>0$ and $M_1,\dots,M_k\in\sM$ such that $U(z,r)\cap M=U(z,r)\cap\bigcup_{i=1}^k M_i$. Since $\widetilde{M}\coloneqq \bigcup_{i=1}^k M_i\in D_2$ by \eqref{sjedno}, and $d_{\widetilde{M}}=d_{M}$ on $U(z,\frac{r}{2})$, we have that $d_{M}$ is DC on $U(z,\frac{r}{2})$. 
	
	Consequently, using also \eqref{loksem}, we obtain that
	$d_M$ is locally DC on $\er^2$, and so it is DC by Lemma~\ref{vldc}~(ii).  
\end{proof}

\section{Every (s)-set is a  $\cal D_2$ set }\label{sd2}

In the next section
we will give a complete characterization of $\cal D_2$ sets using 
``special  $\cal D_2$ sets''
  called {\it (s)-sets}. Their definition is formally rather simple, but it is not easy to prove that
  each (s)-set is a  $\cal D_2$ set. In this section we will prove this important fact
	  using the method of   the proof of  \eqref{hlpz1} from \cite{PZ1} together with  \eqref{hlpz1} and several additional ideas.

 \begin{definition}\label{Ch}
	Let $\varnothing \neq S \subset \R^2$ be a closed set. We say that $S$ is an (s)-set if
 there exists $r>0$ such that
	\begin{equation}\label{Ch1}
	\text{$S \subset \bigcup_{i=1}^{k} \graph f_i$ for some DCR functions $f_i: [0,r] \to \R$
	 with  $(f_i)'_+(0)=f_i(0)=0$}
	\end{equation}
	and 
	\begin{equation}\label{Ch2}
	\text{ $S= \bigcup_{h\in H}\graph h$ for a family $H$ of continuous functions  on $[0,r]$.}
	\end{equation}
\end{definition}

\begin{rem}\label{rem:sSetBasicFacts}
	We will prove some non-trivial properties of (s)-sets in 
	Section 5,  here note only that each (s)-set is clearly nowhere dense and path connected. 
\end{rem}

\begin{rem}\label{rem:changingRforSset}
	Let $S$ be an (s)-set with corresponding $r>0$, functions $f_1,\dots,f_k$ and system $H$ as in Definition~\ref{Ch}. Pick $0<\rho<r$ and put  $\tilde S \coloneqq S\cap ([0,\rho]\times\er)$.  
	 Then clearly $\tilde S$ is an (s)-set (with corresponding functions $\tilde f_i=f_i|_{[0,\rho]}$, $i=1,\dots,k$, and system $\tilde H=\{h|_{[0,\rho]}:h\in H \}$).
\end{rem}

 One from our technical tools is the following easy fact (see \cite[Lemma 3.2]{PZ1}).
\begin{lem}\label{kofu}
	Let $V$ be a closed angle in $\R^2$ with vertex $v$ and measure $0<\alpha< \pi$.
	Then  there exist an affine function $S$ on $\R^2$ and  a concave function $\psi$ on $\R^2$ which is Lipschitz with constant
	$\sqrt{2}\tan (\alpha/2)$ such that  $|z-v| + \psi(z)= S(z),\ z \in V$.
\end{lem}
 We will also use the following ``concave mixing lemma'' (\cite[Lemma 3.1]{PZ1}).
\begin{lem}\label{comix}
	Let $U \subset \R^d$ be an open convex set and let $\gamma: U \to \R$  have finite one-sided directional
	derivatives $\gamma_+'(x,v)$, ($x \in U, \ v\in \R^d$).  Suppose that 
	\begin{equation}\label{nezlo}
	\gamma_+'(x,v) + \gamma_+'(x,-v) \leq 0,\ \ x \in U, v\in \R^d,
	\end{equation}
	and that
	\begin{multline}\label{pokrc}
	\text{ $\graph \gamma$ is covered by graphs of a finite number}\\ 
	\text{of concave functions defined  on $U$.}
	\end{multline}
	Then $\gamma$ is a concave function.
\end{lem}
The core of the present section is the proof of the following lemma which
 easily implies that (s)-sets are $\cal D_2$ sets.
	
\begin{lem}\label{L:main}
	Let $f_1,\dots, f_k$ be DC functions on $\er$   such that each $f_i$ is constant on both $(-\infty, 0]$ and $[1,\infty)$. Let $\varnothing \neq M\subset \R^2$ be a closed set
	and $H$ a system of continuous functions on $\R$ such that
	\begin{equation}\label{eq:inclusion}
	M=\bigcup_{h\in H} \graph h \subset \bigcup_{i=1}^k\graph f_i.
	\end{equation}
	Then $M \in \cal D_2$.
\end{lem}
\begin{proof}
	First observe that each $h\in H$ is constant on both $(-\infty, 0]$ and $[1,\infty)$ by continuity of $h$ and \eqref{eq:inclusion}.
	 
	We will proceed in two steps. 
	
	I)\ \ In the first step we will prove that 
	\begin{multline}\label{gamma}
	\text{there exists a
		concave function $\Gamma$ on $\er^2$}\\
	\text{such that the function $d_M + \Gamma$ is locally
		concave on  $\R^2 \setminus M$.}
	\end{multline}
	
	Observe that Lemma \ref{vldc} (iii) implies that there exists
	 $K>0$ such that the functions $f_1,\dots,f_k$ are $K$-Lipschitz
	 on $[0,1]$ and consequently
	\begin{equation}\label{hkl}
	\text{each function $h\in H$ is $K$-Lipschitz on $[0,1]$}
	\end{equation}
	by \eqref{eq:inclusion} and Lemma \ref{lipmix}.

	If $h\in H$ and $n\in\en$, denote by 
	 $h_n$  the function on $\R$ for which 
	$h_n\left(\frac{i}{n}\right)=h\left(\frac{i}{n}\right)$, $i=0,\dots,n$, which is affine
	on each interval $[\frac{i-1}{n},\frac{i}{n}]$, $i=1,\dots,n$, and which
	is constant on both $(-\infty,0]$ and   $[1,\infty)$.

	For $n\in\en$, 
	 set $H_n\coloneqq \{h_n:\ h\in H\}$ and 
	\begin{equation}\label{mnsj}
	M_n\coloneqq\bigcup_{h\in H_n} \graph h.
	\end{equation}
	 Using   \eqref{eq:inclusion}, we obtain that each $H_n$ is finite and consequently
	 each $M_n$ is closed.
	
	Obviously  $M_n\cap((-\infty,0]\times\er)=M\cap ((-\infty,0]\times\er)$ and 
	$M_n\cap ([1,\infty)
	\times\er)=M\cap ([1,\infty)\times\er)$, $n\in\en$, and
	 \eqref{hkl} easily implies that
	 $M_n\cap ([0,1]\times\er)\to M\cap ([0,1]\times\er)   $ in the Hausdorff metric. Consequenly, we easily obtain that
	$$  d_{M_n} \to d_M\ \ \ \text{on}\  \ \ \R^2.$$

	Now we will show that, to obtain \eqref{gamma}, it is sufficient to find $D>0$ and concave functions
	$\Psi_n$, $n\in \en$, such that 
	\begin{equation}\label{eq:sumOfAngelsFinite}
	\Psi_n \ \ \text{is}\ \ D\text{-Lipschitz on}\ \  \er^2\ \ \text{for each}\ \  n\in \en
	\end{equation}
	and
	\begin{equation}\label{mnpn}
	\text{ $d_{M_n}+\Psi_n$ is locally concave on $\er^2\setminus M_n$ for each $n\in \en$.}
	\end{equation}
	So suppose that such $D$ and $\{\Psi_n\}$, $n\in \en$,  are given
	and consider an arbitrary $z\in \er^2\setminus M$.
	Then there exist $r>0$ and $n_0\in \en$
	such that, for $n \geq n_0$, $B(z,r) \cap M_n = \varnothing$ and, consequently, \eqref{mnpn} easily implies that $d_{M_n}+\Psi_n$ is concave on $B(z,r)$. 
	
	Further observe that we can suppose  that $\Psi_n(0)=0$, $n\in \en$. Then \eqref{eq:sumOfAngelsFinite} gives that the sequence $\{\Psi_n\}$ is equicontinuous and pointwise bounded, and  we can use a well-known version of Arzel\`a-Ascoli theorem (see e.g. \cite[Theorem 4.44, p. 137]{Fo}) to obtain a subsequence $\{\Psi_{n_p}\}$ converging to a $D$-Lipschitz concave function $\Gamma$ on $\er^2$.
	Then $d_{M_{n_p}}+\Psi_{n_p}\to d_M+\Gamma$, consequently
	$d_M+\Gamma$ is  concave on $B(z,r)$, and so \eqref{gamma} holds.
	
	To find $\{\Psi_n\}$ and $D$, we first
	define a finite subset $A_n$ of $M_n$ by
	$$ A_n\coloneqq  \left\{ \left(\frac{i}{n}, h\left(\frac{i}{n}\right)\right):\ 0\leq i \leq n,\ h \in H\right\}.$$
	For each $n\in \en$ and $a=(a_1,a_2) \in A_n$ set
	$$ s^+(n,a)\coloneqq  \max \{h'_+(a_1):\ h\in H_n, h(a_1)=a_2\},\ s_+(n,a)\coloneqq  \min \{h'_+(a_1):\ h\in H_n, h(a_1)=a_2\},$$
	$$ s^-(n,a)\coloneqq  \max \{h'_-(a_1):\ h\in H_n, h(a_1)=a_2\},\ s_-(n,a)\coloneqq  \min \{h'_-(a_1):\ h\in H_n, h(a_1)=a_2\}.$$

	% \begin{equation*}
	%\max \{s^+(n,a) - s_-(n,a), s^-(n,a) - s_+(n,a)\}  \leq \frac{f(a_1 + 1/n)- f(a_1)}{1/n}
	%- \frac{f(a_1)- f(a_1-1/n)}{1/n}
	%\end{equation*}	
	Now we will prove that there is a constant $C>0$  such that, for each $n \in \en$,
	\begin{equation}\label{eq:sumOf2ndDiffFinite}
	\sum_{a \in A_n} |s^+(n,a) - s_-(n,a)| \leq C,\ \sum_{a \in A_n}  |s^-(n,a) - s_+(n,a)| \leq C.
	\end{equation}
	To this end, consider a convex function $\vf$ which corresponds to 
	    $J\coloneqq \R$ and $F_i\coloneqq f_i$, $i=1,\dots,k$,   
	 by Lemma~\ref{vmix}.
	Choose $L>0$ such that $\vf$ is $L$-Lipschitz on $[-1,2]$.
	Further consider arbitrary $n \in \en$, $0\leq i \leq n$ and $a= (a_1,a_2) \in A_n$
	with $a_1=i/n$. 
	Now choose $\tilde h \in H$ and $\hat h \in H$ such that
	$$ \tilde h(i/n)= \hat h (i/n) = a_2,\ \ (\tilde h_n)_+'(i/n)= s^+(n,a),\ \ (\hat h_n)'_-(i/n)= s_-(n,a).$$
	Set  $g(x)\coloneqq \tilde h(x)$ for $x \geq i/n$ and  $g(x)\coloneqq \hat h(x)$ for $x < i/n$. 
	Then clearly
	$$  |s^+(n,a) - s_-(n,a)| = \left| \frac{g\left( \frac{i+1}{n}\right) -g\left(\frac{i}{n}
		\right)}
	{\frac{1}{n}} -  \frac{g\left( \frac{i}{n} \right) -g\left(\frac{i-1}{n} 
		\right)}
	{\frac{1}{n}}         \right|.    $$ 
	Since $\graph g \subset \bigcup_{i=1}^k\graph f_i$, by its choice, $\vf$ is a control function for $g$
	and so Lemma~\ref{smercontr} and the above equality imply 
	\begin{equation}\label{eq:2L}
	|s^+(n,a) - s_-(n,a)| \leq  
	\frac{\vf\left( \frac{i+1}{n}\right) -\vf\left(\frac{i}{n}
		\right)}
	{\frac{1}{n}} -  \frac{\vf\left( \frac{i}{n} \right) -\vf\left(\frac{i-1}{n} 
		\right)}
	{\frac{1}{n}}\leq 2L.
	\end{equation}
	Consequently
	\begin{multline*}
	\sum_{a \in A_n} |s^+(n,a) - s_-(n,a)| \leq \sum_{i=0}^n \ \ \sum_{(a_1,a_2) \in A_n, a_1=\frac{i}{n}}
	\left( \frac{\vf\left( \frac{i+1}{n}\right) -\vf\left(\frac{i}{n}
		\right)}
	{\frac{1}{n}} -  \frac{\vf\left( \frac{i}{n} \right) -\vf\left(\frac{i-1}{n} 
		\right)}
	{\frac{1}{n}}         \right)\\
	\leq k \sum_{i=0}^n 
	\left( \frac{\vf\left( \frac{i+1}{n}\right) -\vf\left(\frac{i}{n}
		\right)}
	{\frac{1}{n}} -  \frac{\vf\left( \frac{i}{n} \right) -\vf\left(\frac{i-1}{n} 
		\right)}
	{\frac{1}{n}}         \right)
	\\
	= k  \left( \frac{\vf\left( \frac{n+1}{n}\right) -\vf\left(1
		\right)}
	{\frac{1}{n}} -  \frac{\vf\left(0 \right) -\vf\left(-\frac{1}{n} 
		\right)}
	{\frac{1}{n}}         \right) \leq 2Lk  \eqqcolon C.
	\end{multline*}
	The second inequality of    \eqref{eq:sumOf2ndDiffFinite} follows quite analogously.
	
	   For each $n \in \en$ and $a=(a_1,a_2) \in A_n$, set   
	$$ p^+(n,a)\coloneqq  (a_1+ 1/n, a_2 + s^+(n,a)/n),\ p_+(n,a)\coloneqq  (a_1+ 1/n, a_2 + s_+(n,a)/n),$$
	$$  p^-(n,a)\coloneqq  (a_1- 1/n, a_2 - s^-(n,a)/n),\ p_-(n,a)\coloneqq  (a_1- 1/n, a_2 - s_-(n,a)/n),$$
	$$A_n^1\coloneqq  \{a\in A_n:\ s^+(n,a) - s_-(n,a)>0\},\ \  A_n^2\coloneqq  \{a\in A_n:\ s^-(n,a) - s_+(n,a)>0\}.$$   
	Further set
	\begin{equation*}
	V^1_{n,a}\coloneqq\{z\in\er^2: 
	\langle z - a ,p^+(n,a)- a\rangle\leq 0,\;
	\langle z - a, p_-(n,a) -a \rangle\leq 0 
	\}\ \ \text{if}\ \ a \in A_n^1
	\end{equation*}
	and
	\begin{equation*}
	V^2_{n,a}\coloneqq\{z\in\er^2: 
	\langle z - a ,p_+(n,a)- a\rangle\leq 0,\;
	\langle z - a, p^-(n,a) -a \rangle\leq 0 
	\}\ \ \text{if}\ \ a \in A_n^2.
	\end{equation*}
	It is easy to see that each $V^1_{n,a}$ (resp. $V^2_{n,a}$) 
	is a closed angle with vertex $a$ and measure 
	$$
	\alpha^1(n,a)\coloneqq \arctan s^+(n,a) - \arctan s_-(n,a)\in (0,\pi)
	$$
	 
	$$(\text{resp.}\quad\alpha^2(n,a)\coloneqq \arctan s^-(n,a) - \arctan s_+(n,a) \in (0,\pi)).$$
	For  $a \in A_n^1$ (resp. $a \in A_n^2$)  
	let $\psi^1_{n,a}$ and $S^1_{n,a}$ (resp. $\psi^2_{n,a}$ and $S^2_{n,a}$) be the concave and affine functions on
	$\R^2$ which correspond  to $V^1_{n,a}$ (resp. $V^2_{n,a}$) by Lemma \ref{kofu}. 
	If $a \in A_n \setminus A_n^1$ (resp. $a \in A_n \setminus A_n^2$  ),  put $\psi^1_{n,a}(z)\coloneqq 0$ and $S^1_{n,a}(z)\coloneqq 0$ (resp. $\psi^2_{n,a}(z)\coloneqq 0$ and $S^2_{n,a}(z)\coloneqq 0$), $z \in \R^2$. Set
	$$ \Psi_n\coloneqq\sum_{a \in A_n} (\psi^1_{n,a}+ \psi^2_{n,a}).$$
	Now fix an arbitrary $a \in A_n^1$. Using \eqref{eq:2L} we easily obtain 
	\begin{equation*}
	\alpha^1(n,a) \leq s^+(n,a)-s_-(n,a)\leq 2L.
	\end{equation*}
	Further, since the function $\tan$ is convex on $[0,\pi/2)$, the function $\omega(x)= \frac{\tan x}{x}$
	is increasing on $(0,\pi/2)$. These facts easily imply 
	$$   \sqrt{2}\tan\left(\frac{\alpha^1(n,a)}{2}\right)\leq  \sqrt{2}\cdot\frac{\alpha^1(n,a)}{2} \cdot \frac{L}{ \arctan L}
	\leq (s^+(n,a)-s_-(n,a))\cdot \frac{L}{ \sqrt{2}\arctan L}.$$
	So, by the choice of $\psi^1_{n,a}$, we have that
	\begin{equation}\label{odhlk1}
	\psi^1_{n,a}  \ \ \text{is Lipschitz with constant}\ \   (s^+(n,a)-s_-(n,a))\cdot \frac{L}{ \sqrt{2}\arctan L}.
	\end{equation}
	Quite similarly we obtain that, for each $a \in A_n^2$, 
	\begin{equation}\label{odhlk2}
	\psi^2_{n,a}  \ \ \text{is Lipschitz with constant}\ \ (s^-(n,a)-s_+(n,a))\cdot \frac{L}{\sqrt{2} \arctan L}.      
	\end{equation} 
	Consequently  \eqref{odhlk1},  \eqref{odhlk2}  and    \eqref{eq:sumOf2ndDiffFinite} easily imply
	that there is a constant $D>0$ such that
	\eqref{eq:sumOfAngelsFinite}
	holds.

	To prove   \eqref{mnpn},  it is enough to prove that
	\begin{equation}\label{hnpsin}
	\text{$d_{\graph h}+\Psi_n$
		is locally concave on $\er^2\setminus \graph h$ for each $n \in \en$ and $h \in H_n$.}
	\end{equation}
	
	Indeed, by   \eqref{mnsj}   it is easy to see that
	\begin{equation*}
	d_{M_n}+\Psi_n=\min_{h\in H_n}(d_{\graph h}+\Psi_n)
	\end{equation*}
	and it is enough to use  (on each open ball $U \subset \R^2 \setminus M_n$)  the fact that the minimum of a finite system of concave functions is a concave function.
	
	To prove \eqref{hnpsin}, fix an arbitrary $n\in \en$ and  $h \in H_n$.
	
	For  $i=-1,\dots,n+1$ denote $z_i: =\left(\frac{i}{n},h\left(\frac{i}{n}\right)\right)$ and
	\begin{equation*}
	V_i\coloneqq  \{ z \in \R^2:\ \langle z-z_i, z_{i+1}-z_{i}\rangle \leq 0,\ \langle z-z_i, z_{i-1}-z_{i}\rangle \leq 0 \}.
	\end{equation*}
	
	Now, for a fixed $i$, denote $a\coloneqq z_i$. 
	Then clearly $a\in A_n$
	and, if the points  $z_{i-1},z_{i}, z_{i+1}$ are not collinear, then
	\begin{equation}\label{inklv}
	\text{ either $a \in A^1_n$ and $V_i \subset V^1_{n,a}$, or   $a \in A^2_n$  and  $V_i \subset V^2_{n,a}$.}
	\end{equation}
	Indeed, observe that
	$$  s_+(n,a) \leq  \frac{h\left(\frac{i+1}{n}\right) - h\left(\frac{i}{n}\right)}{\frac{1}{n}} \leq
	s^+(n,a),\ \ \  s_-(n,a) \leq  \frac{h\left(\frac{i}{n}\right) - h\left(\frac{i-1}{n}\right)}{\frac{1}{n}} \leq
	s^-(n,a).$$
	So, if $\frac{h\left(\frac{i+1}{n}\right) - h\left(\frac{i}{n}\right)}{\frac{1}{n}} >
	\frac{h\left(\frac{i}{n}\right) - h\left(\frac{i-1}{n}\right)}{\frac{1}{n}}$, then $a \in A^1_n$
	and an easy geometrical observation shows that $V_i \subset V^1_{n,a}$. Similarly, if
	$\frac{h\left(\frac{i+1}{n}\right) - h\left(\frac{i}{n}\right)}{\frac{1}{n}} <
	\frac{h\left(\frac{i}{n}\right) - h\left(\frac{i-1}{n}\right)}{\frac{1}{n}}$, then $a \in A^2_n$
	 and $V_i \subset V^2_{n,a}$.

	Denote 
	$l_{-}\coloneqq \er\times \{h(0)\}$,
	$l_{+}\coloneqq \er\times \{h(1)\}$
	and for $n\in\en$ and $i=0,\dots,n-1$ denote
	\begin{equation*}
	\eta_i\coloneqq d_{l(z_i,z_{i+1})}+\Psi_{n}
	\end{equation*}
	and
	\begin{equation*}
	\eta_{-1}\coloneqq d_{l_-}+\Psi_n,\quad
	\eta_{n}\coloneqq d_{l_+}+\Psi_n.
	\end{equation*}
	We will prove that $\gamma\coloneqq  d_{\graph h}+\Psi_n$
	is locally concave on $\er^2\setminus \graph h$ using Lemma \ref{comix}.
	So fix an arbitrary $b\in \er^2\setminus \graph h$ and $\delta>0$ such that $U\coloneqq  U(b, \delta)
	\subset \er^2\setminus \graph h$. Then condition \eqref{nezlo} of Lemma \ref{comix} holds by
	\eqref{zlose} and \eqref{loksem}. To prove condition \eqref{pokrc}, consider an arbitrary $z\in U$ and choose
	$z^*\in \graph h$ such that $d_{\graph h}(z)=|z-z^*|$.
	
	First note that if $z^*\not\in\{z_{0},\dots,z_{n}\}$ or $z^* = z_i$ for some
	$i \in \{0,1,\dots,n\}$ and the points    $z_{i-1},z_{i}, z_{i+1}$ are collinear,
	then clearly
	\begin{equation}\label{ety}
	\gamma(z)= d_{\graph h}(z)+\Psi_n(z)\in 
	\bigcup_{i\in\{-1,\dots,n\}} \{\eta_{i}(z)\}.
	\end{equation}
	Further assume that $z^*=z_i$ for some $i\in\{0,\dots,n\}$ and the points  $z_{i-1},z_{i}, z_{i+1}$ are not  collinear. 
	Then clearly $z\in V_i$ and  $\eqref{inklv}$ holds. Consequently
	$$ \text{either \ \ $|z-z^*| + \psi^1_{n, z_i} (z) = S^1_{n, z_i} (z)$\ \  or \ \  
		$|z-z^*| + \psi^2_{n, z_i} (z) = S^2_{n, z_i} (z)$},$$
	and therefore
	\begin{equation}\label{kdezh}
	\gamma(z)= d_{\graph h}(z) + \Psi_n(z) \in \{ S^1_{n, z_i} (z) + (\Psi_n -\psi^1_{n, z_i}) (z) ,
	S^2_{n, z_i} (z) + (\Psi_n -\psi^2_{n, z_i})(z)\}.
	\end{equation}
	Since the graph of each function
	$\eta_{i}$, $-1\leq i \leq n$, can be clearly covered by graphs of two concave functions
	and the functions $\Psi_n -\psi^1_{n, z_i}, \Psi_n -\psi^2_{n, z_i}$  ($i=0,\dots,n$)
	   are  concave,
	\eqref{ety}     and       \eqref{kdezh}  imply  \eqref{pokrc} and so $\gamma$ is concave on $U$ and therefore \eqref{hnpsin} holds, which completes the proof   of  \eqref{gamma}. 
	\medskip
	
	II)\ \ In the second step we first observe that
	by  \eqref{hlpz1}  there exist concave functions $\omega_i$, $1\leq i \leq k$,  such
	that each function $d_{\graph f_i} + \omega_i$ is concave  on $\R^2$.  Set  
	$$  \omega\coloneqq  \sum_{i=1}^k \omega_i\ \ \text{and}\ \ \sigma\coloneqq   \Gamma + \omega.$$
	Then
	\begin{equation}\label{zopak} 
	d_{\graph f_i} + \sigma  \ \ \text{is concave},\ \  1\leq i \leq k
	\end{equation}
	and by \eqref{gamma}
	\begin{equation}\label{zopak2}
	d_M + \sigma\ \ \text{is locally concave on}\ \ \R^2\setminus M.
	\end{equation}

	It is sufficient to prove that $d_M +\sigma$ is concave on $\R^2$.

	For $i,j\in\{1,\dots,k\}$, $i\not=j$, denote by $P_{i,j}$ the set of all $x\in\R$ such that $f_i(x)=f_j(x)$ and such that for every $\eps>0$ there is $z\in(x-\eps,x+\eps)$ satisfying $f_i(z)\not=f_j(z)$. Obviously, each  $P_{i,j}$ is a closed nowhere dense set and  $P_{i,j}\subset [0,1]$.
	
	Set
	$$  P_i\coloneqq  \bigcup \{P_{i,j}:\ 1\leq j \leq k,\ j\neq i\},\ \ P\coloneqq  \bigcup_{i=1}^k  P_i,$$
	$$  P^*_i\coloneqq   \{(x, f_i(x)):\ x \in  P_i\},\ \  P^*\coloneqq  \bigcup_{i=1}^k  P^*_i.$$

	Note that for every $z\in M\setminus P^*$ there is $i\in\{1,\dots,k\}$ and $\rho>0$ such that $M\cap U(z,\rho)=\graph f_i\cap U(z,\rho)$ and so
	\begin{equation}\label{rednaf}
	d_M(u) = d_{\graph f_i}(u),\ \ u \in  U(z,\rho/2).
	\end{equation}
	
	To prove the concavity of $d_M + \sigma$, it is clearly sufficient to prove that for each $p,q \in \R^2$
	and $\ep>0$ there exists a line $l$ which meets both $U(p,\ep)$ and $U(q,\ep)$ and
	$d_M +\sigma$ is concave on $l$. To this end, choose arbitrary $p,q, \ep$. Further, using the notation
	$$l(m,c)\coloneqq   \{(x,y):\ y= mx+c\},\ \ m,c \in \R,$$
	we  choose a line $l(m_0,c_0)$ which meets both $U(p,\ep)$ and $U(q,\ep)$. 
	
	Now observe that for each $c \in \R$ we have $l(m_0, c) \cap P_i^* \neq \varnothing$ if and only if
	there is $x \in P_i$ such that $m_0x +c= f_i(x)$, i.e. $c \in g_i(P_i)$, where
	$g_i(x)\coloneqq  f_i(x)-m_0x, x \in \R$.  Since $g_i$ is DC and $P_i \subset [0,1]$ is nowhere dense, Lemma \ref{obrdc}
	implies that $C_i \coloneqq  \{c\in \R:\ l(m_0,c) \cap P_i^* \neq \varnothing\}$ is nowhere dense.
	Consequently we can choose $c_1 \in \R$ such that $l\coloneqq  l(m_0,c_1) \subset \R^2 \setminus P^*$
	and $l$  meets both $U(p,\ep)$ and $U(q,\ep)$. 
	By    \eqref{rednaf}, \eqref{zopak} and  \eqref{zopak2} we obtain that $d_M + \sigma$ is locally concave at each point
	of $l$, and thus concave on $l$.
\end{proof}

%\begin{cor}\label{cor:emptyInterior}
%	Let $M$ be a (s)-set, then $d_M$ is DC.
%\end{cor}

\begin{cor}\label{cor:CHsetDistanceDC}
	Let $M\subset\er^2$ be as in Lemma~\ref{L:main}. 
	Then $\widetilde M\coloneqq M\cap([0,1]\times\er)\in\D_2$.
\end{cor}

\begin{proof}
	Let $f_1,\dots,f_k$ and $H$ be as in Lemma~\ref{L:main}.
	First note that, by Lemma~\ref{L:main} and \eqref{loksem}, $d_{\widetilde M}$ is locally DC on $$\er^2\setminus (( M\cap(\{0\}\times\er))\cup (M\cap(\{1\}\times\er)))
	\eqqcolon\er^2\setminus(M_0\cup M_1). $$
	
	We prove that $d_{\widetilde M}$ is DC on some neighbourhood of each point in  $M_0$.	
	To do that pick some $z\in M_0$.
	Let $H_z$ be the system of all functions $h\in H$ such that $(0,h(0))=z$ and put
	\begin{equation*}
	f_z(x)=\max_{h\in H_z} h(x),\quad g_z(x)=\min_{h\in H_z} h(x), \quad x\in [0,\infty),
	\end{equation*}
	and 
	$$
	D_z\coloneqq \{(x,y):x\in[0,\infty),\; g_z(x)\leq y\leq f_z(x) \}.
	$$
	
	Note that $d_{M}$ is DC by Lemma~\ref{L:main}.
	  Since the functions $f_1,\dots,f_k$ are Lipschitz on $[0,1]$ by   Lemma \ref{vldc} (iii), they are Lipschitz on $[0, \infty)$. Consequently both $f_z$ and $g_z$ are Lipschitz by  \eqref{eq:inclusion} and  Lemma~\ref{lipmix},     and therefore they are DCR by \eqref{eq:inclusion}  and Lemma~\ref{L:DCRmixing}.  
 Therefore $d_{D_z}$ is DC by 
	\eqref{eq:boundaryEquivalenceGeneral}   and \eqref{loksjgr}, since
	 $\partial D_z$ is clearly the union of a locally finite system
	 of DC graphs. 
	It is easy to see (using the fact that the set $M_0$ has cardinality at most $k$ and so, in particular, is finite) that for sufficiently small $\rho>0$ we have $d_{\widetilde M}(w)\in\{d_{D_z}(w),d_M(w)\}$ whenever $w\in U(z,\rho)$ and so $d_M$ is DC on $U(z,\rho)$ by Lemma~\ref{vldc}~(vi).
	
	Quite analogically one can prove that $d_{\widetilde M}$ is also DC on a neighbourhood of each point of $M_1$ and so $d_{\widetilde{M}}$ is locally DC and therefore (by Lemma~\ref{vldc}~(ii)) DC on $\er^2$ and $\widetilde{M}\in\D_2$.
\end{proof}

\begin{cor}\label{cor:sSetInD2}
	Every (s)-set $S\subset \R^2$ is a $\cal D_2$ set.
\end{cor}
\begin{proof}
	Let $S$ be an (s)-set and let $r>0$, $f_1,\dots,f_k$ and $H$ be as in Definition~\ref{Ch}.
	We may assume  (applying a suitable similarity and using Remark 
 \ref{rem:satbilityOfD2} (i) if necessary)   that $r=1$.
	For $h\in H$ define $\tilde{h}:\er\to\er$ by $\tilde h=h$ on $[0,1]$, $\tilde h=h(0)$ on $(-\infty,0]$ and $\tilde{h}=h(1)$ on $[1,\infty)$. 
	Set $\tilde H: = \{\tilde h:\ h\in H\}$.
	Similarly we extend functions $f_i$ calling the extensions $\tilde f_i$.
	Clearly (by Lemma~\ref{vldc}~(vi)) each $\tilde f_i$ is a DC function on $\er$.
	Put $M\coloneqq\bigcup_{h\in \tilde H}\graph h$.
	Then
	\begin{equation*}
	M=\bigcup_{h\in \tilde H} \graph h \subset \bigcup_{i=1}^k\graph \tilde f_i.
	\end{equation*}
 	Since $S=M\cap([0,1]\times\er)$ and  $M$ satisfies  the assumptions of Lemma~\ref{L:main}, we obtain $S\in\D_2$ by Corollary~\ref{cor:CHsetDistanceDC}. 
\end{proof}

\section{Complete characterizations of $\cal D_2$ sets}
%Characterization of $\cal D_2$ sets by (s)-sets and by nowhere dense $\cal D_2$ sets}

\begin{lem}\label{L:sequenceToFunction}
	Let $\{a_n\}_{n=1}^\infty$ and $\{A_n\}_{n=1}^\infty$ be sequences of real numbers such that
	\begin{enumerate}
		\item[(i)] $0<a_{n+1}\leq\frac{a_n}{3}$, $n\in\en$, and
		\item[(ii)] $\sum_{n=1}^{\infty}\frac{|A_n|}{a_n}<\infty$.
	\end{enumerate}
	Then the function 
	\begin{equation*}
	f(x)=\begin{cases}
	\frac{A_n-A_{n+1}}{a_n-a_{n+1}}(x-a_{n+1}) + A_{n+1} &\quad\text{if $x\in(a_{n+1},a_n]$}\\
	0&\quad \text{if $x= 0$}
	\end{cases}
	\end{equation*}
	is a DCR function on $[0,a_1]$.
\end{lem}

\begin{proof}
	First note that condition (i) implies
	\begin{equation}\label{eq:inequalitiesa_n}
	a_n-a_{n+1}\geq a_n-\frac{a_n}{3} = \frac23 a_n\geq 2 a_{n+1}.
	\end{equation}
	Since $a_n \to 0$ by (i), we obtain $A_n = f(a_n) \to 0$ by (ii) and so
	$f$ is continuous. 
	Clearly $f'_-(x) = f'_-(a_n) = \frac{A_n-A_{n+1}}{a_n-a_{n+1}}$ for $x \in (a_{n+1},a_n]$, $n \in \en$.
	Using (ii) and \eqref{eq:inequalitiesa_n}, we obtain
	\begin{equation*}\label{sumder}
	\sum_{n=1}^{\infty}|f_{-}'(a_n)|\leq \sum_{n=1}^{\infty}\left(\left|\frac{A_n}{a_n-a_{n+1}}\right|+ \left|\frac{A_{n+1}}{a_n-a_{n+1}}\right|\right)
	\leq \sum_{n=1}^{\infty}\left(\frac{|A_n|}{\frac23a_n}+ \frac{|A_{n+1}|}{2a_{n+1}}\right) < \infty.
	\end{equation*}
	Therefore we easily obtain
	$$ V(f'_-,(0,a_1))= \sum_{n=1}^{\infty} |f'_-(a_n) - f'_-(a_{n+1})| < \infty,$$
	and so $f$ is a DCR function by  Lemma  \ref{dcr}.
\end{proof}

For $r>0$ put $$A_{r}^u=\{(x,y): 0\leq x\leq r,\;  |y| \leq ux \},\quad A^u=\{(x,y): 0\leq x,\;  |y| \leq ux \},$$
 $$S_{r}^u=\{(x,y): |x|\leq r,\;  |y| \leq u|x| \}\quad\text{and}\quad S^u=\{(x,y): x\in\er,\;  |y| \leq u|x| \}.$$

In the proof    of  Lemma \ref{uhldiry}    we will use the following geometrically obvious lemma.

\begin{lem}\label{geom}
	Let $u>0$. Then there exists $\alpha>0$ such that 
	$ d_M(x,y) \geq \alpha \cdot \rho$, whenever $M \subset \R^2$, $z=(a,b) \in \R^2$, $\rho>0$ and either
	\begin{equation}\label{plu}
	M \cap (z + A_{\rho}^{3u}) \subset \{z\},\ (x,y) \in (z+A_{\rho}^{2u}),\ x = a + \frac \rho 2,
	\end{equation}
	or
	\begin{equation}\label{min}
	M \cap (z - A_{\rho}^{3u}) \subset \{z\},\ (x,y) \in (z -A_{\rho}^{2u}),\ x = a - \frac \rho 2.
	\end{equation}
\end{lem}

\begin{lem}\label{uhldiry}
	Let $M \subset \R^2$ be closed, $z \in M$ and $u>0$. Let there exist sequences $\{z_n\}$ in $M$ and $\{\rho_n\}$ in $(0,\infty)$
	such that $z_n \to z$  and, for each $n \in \en$,
	\begin{enumerate}
		\item[(i)]
		$z_n \in  (z + S^u)\setminus \{z\}$ and
		\item[(ii)]
		either  $(z_n + A_{\rho_n}^{3u}) \cap M = \{z_n\}$\ or\ $(z_n - A_{\rho_n}^{3u}) \cap M = \{z_n\}$.
	\end{enumerate}
	Then  $M \notin \cal D_2$.
\end{lem}
\begin{proof}
	Suppose, to the contrary, that $d_M$ is DC.
	
	We can suppose $z=0$. Let $z_n = (a_n,b_n)$. 
	Without  any loss of generality we can suppose $a_1> a_2>\dots>0$. To see this, we can pass to a subsequence and work
	with $M^s\coloneqq  \{(x,y): (-x,y) \in M\}$ (which belong to  $\cal D_2$
	if and 
	only if $M \in \cal D_2$) instead of $M$, if necessary.
	Further, we can assume (passing several times to a subsequence), that
	\begin{equation}\label{prpod}
	a_{n+1}\leq \frac{a_n}{3},\ \ \text{for each}\ n \in \en,
	\end{equation}
	and, for some $K\in [-u,u]$,
	\begin{equation}\label{drpod}
	\frac{b_n}{a_n} \to K\ \ \ \text{and}\ \ \ \sum_{n=1}^{\infty} \left| \frac{b_n}{a_n} - K\right| \ < \ \infty.
	\end{equation}
	Note that also (by $z_n\in A^u$ and \eqref{prpod})
	\begin{equation}\label{prpod2}
\left| \frac{b_n- b_{n+1}}{a_n-a_{n+1}}\right|
\leq \frac{a_n u +\frac{1}{3}a_{n}u}{\frac23a_n}
\leq 2u 
\ \ \ \text{for each} \ n \in \en.
	\end{equation}

	Set  $A_n\coloneqq  b_n-a_nK$. Then  assumptions  (i) and (ii) of Lemma~\ref{L:sequenceToFunction} are satisfied by \eqref{prpod} and \eqref{drpod}, and consequently we know that
	the function
	\begin{equation*}
	f(x)=\begin{cases}
	\frac{A_n-A_{n+1}}{a_n-a_{n+1}}(x-a_{n+1}) + A_{n+1} &\quad\text{if $x\in(a_{n+1},a_n]$}\\
	0&\quad \text{if $x= 0$}
	\end{cases}
	\end{equation*} 
	is a DCR function on $[0,a_1]$ and  Lemma \ref{vldc} (i),(ii),(iv)  easily imply that the functions $g(x)=f(x)+Kx,\ x \in [0,a_1]$, and $F(x)=d_M(x,g(x)),\  x \in [0,a_1]$, are also DCR functions. Observe that 
	\begin{equation}\label{vlg}
	\text{	$g(a_n)= b_n,\ n\in \en$ and $g$ is linear on each interval
		$[a_{n+1},a_n]$.}
	\end{equation}
	Further, $F(0)=0$ and $F(a_n)=0$, $n \in \en$. So Lemma \ref{strict} implies that
	\begin{equation}\label{strnul}
	\text{$0$ is the right strict derivative of $F$ at $0$.}
	\end{equation}
	Now consider $n>1$ and choose an $\alpha >0$ which corresponds to our $u$ by Lemma \ref{geom}. 
	
	If  $(z_n + A_{\rho_n}^{3u}) \cap M = \{z_n\}$, choose $0<r_n< \rho_n$ such that
	$a_n+ r_n < a_{n-1}$ and set $x_n\coloneqq  a_n + r_n/2$, $y_n\coloneqq  g(x_n)$.  Since \eqref{vlg} and  \eqref{prpod2}   imply 
	$(x_n, y_n) \in  z_n + A_{r_n}^{2u}$, we can apply Lemma \ref{geom} (with $z=0$, $\rho=r_n$, $x=x_n$, $y=y_n$) and
	obtain  $d_M(x_n,y_n) = F(x_n) \geq \alpha r_n$, Consequently 
	\begin{equation}\label{velsmer}
	\frac{F(x_n)-F(a_n)}{x_n-a_n} \geq 2 \alpha.
	\end{equation}
	
	If $(z_n - A_{\rho_n}^{3u}) \cap M = \{z_n\}$,  choose $0<r_n< \rho_n$ such that
	$a_n- r_n > a_{n+1}$ and set   $x_n\coloneqq  a_n - r_n/2$,   $y_n\coloneqq  g(x_n)$. In the same way as in the first case we also
	obtain \eqref{velsmer}. 
	
	Now observe that, by the definition of the strict right derivative, \eqref{velsmer} contradicts \eqref{strnul}.
	
\end{proof}

\begin{lem}\label{L:coneProjection}
	Let $M\in \cal D_2$, $z=(x,y) \in M$, $s>0$ and $u>0$.
	Then the following assertions hold.
	\begin{enumerate}
		\item  [(i)]
		If 
		\begin{equation}\label{prava}
		\partial M\cap (z+ A_s^{3u})
		\subset z+ A_s^u,
		\end{equation}\label{leva}
		then there is $s\geq r>0$ such that either $M\cap (z+A_r^{3u})=\{z\}$, or 
		$\pi_1(M\cap (z+A_r^{3u})    )=[x,x+r]$.
		\item  [(ii)]
		If 
		\begin{equation}
		\partial M\cap (z- A_s^{3u})
		\subset z- A_s^u,
		\end{equation}
		then there is $s\geq r>0$ such that either $M\cap (z-A_r^{3u})=\{z\}$, or 
		$\pi_1(M\cap (z-A_r^{3u})    )=[x-r,x]$. 
	\end{enumerate}
\end{lem}
\begin{proof}
	We will prove only assertion (i); the proof of (ii) is quite analogous.
	
	Set  $K\coloneqq  \pi_1(M\cap (z+A_s^{3u})    )$ and observe that condition \eqref{prava}
	implies that either $\pi_1(M\cap (z+A_s^{3u})    )=[x,x+s]$  or
	\begin{equation}\label{kdem}
	M\cap (z+ A_s^{3u})
	\subset z+ A_s^u.
	\end{equation}
	So we can suppose that \eqref{kdem} holds.
	
	Now suppose to the conratry that no $s\geq r >0$ from the assertion of the lemma exists.
	Since $K$ is compact, we can clearly find sequences of positive numbers $\{x_n\}$, $\{\rho_n\}$
	such that $x_n \to 0$,
	\begin{equation}\label{vlxrho}
	x_n \in K\ \ \ \text{and}\ \ \ (x_n,x_n+\rho_n) \cap K = \varnothing\ \ \ \text{for each}\ \ \ n\in \en.
	\end{equation}
	By the definition of $K$  and \eqref{kdem} there exist $y_n$, $n\in \en$, such that 
	$z_n\coloneqq  (x_n,y_n) \in M \cap (z + A_s^u)$. Since \eqref{vlxrho} clearly implies
	$(z_n + A_{\rho_n}^{3u}) \cap M = \{z_n\}$, we obtain by Lemma \ref{uhldiry} that $M \notin \cal D_2$,
	which is a contradiction.
\end{proof}

\begin{lem}\label{L:localExtension}
	Let $M\in \cal D_2$  and $s,u>0$.
	Suppose that $0\in M$ and 
	$f_1,\dots, f_k:[0,s]\to\er$ are  $u$-Lipschitz functions
	such that $f_i(0)=0$, $i=1,\dots,k$, and
	\begin{equation}\label{eq:twoInclusionsOne}
	\partial M\cap A_s^{3u}\subset\bigcup_{i=1}^k  \graph f_i.
	\end{equation}
	Let $0$ be an accumulation point of $\partial M \cap A_s^u$.
	Then there is some $0<\rho<s$ such that  for every 
	$$(x,y)\in M \cap \bigcup_{i=1}^k  \graph f_i$$
	with $x\in (0,\rho)$
	there exist $\delta>0$ and a $u$-Lipschitz function $g:[x-\delta,x+\delta]\to\er$ such that $g(x)=y$ and $\graph g \subset M\cap \bigcup_{i=1}^k \graph f_i$.
\end{lem}

\begin{proof}
	%	Pick $s\geq r>0$ as in Lemma~\ref{L:coneProjection}.
	For the sake of brevity, we set $M^*\coloneqq  M \cap \bigcup_{i=1}^k  \graph f_i$ and observe that
	the properties of $f_i$ imply $M^* \subset A_s^u$.
	
	Now consider an arbitrary $z=(x,y) \in M^*$ with $x \in (0,s)$.
	Since $z \in A_s^u$, it is easy to see that we can assign to $z$ a number $R_z>0$
	such that 
	\begin{equation}\label{uhvuh}
	z+ S_{R_z}^{3u} \subset A_s^{3u}
	\end{equation}
	and
	$$  \graph f_i \cap (z+ S_{R_z}^{3u}) = \varnothing \ \ \text{whenever}\ \ 1\leq i \leq k\ \text{and}\ f_i(x) \neq y.$$
	Then, using $u$-Lipschitzness of all $f_i$ and   \eqref{eq:twoInclusionsOne}, we
	obtain
	\begin{equation}\label{eq:insideCone}
	\partial M\cap (z+ S^{3u}_{R_z})\subset z+S^u_{R_z}. 
	\end{equation}
	So, by
	Lemma~\ref{L:coneProjection}, we can choose $0<r_z \leq R_z$ such that
	$\pi_1(M\cap (z+ S^{3u}_{r_z}))$ is one of the following sets:
	$$ \{x\},\ [x-r_z,x],\ [x,x+r_z],\ [x-r_z,x+r_z].$$
	Using Lemma \ref{uhldiry}, we easily obtain that there exists
	$s\geq \rho>0$ such that 
	\begin{equation}\label{eq:fullProjection}
	\pi_1(M\cap (z+ S^{3u}_{r_z}))=[x-r_z,x+r_z],\quad\text{whenever}\quad x\in(0,\rho).
	\end{equation}
	We claim that even
	\begin{equation}\label{eq:fullProjectionWifhFunctions}
	\pi_1(M^*\cap (z+ S^{3u}_{r_z}))=[x-r_z,x+r_z],\quad\text{whenever}\quad x\in(0,\rho).
	\end{equation}
	Indeed, pick $t\in [x-r_z,x+r_z]$. To prove $t \in \pi_1(M^*\cap (z+ S^{3u}_{r_z}))$, we distinguish two cases.
	If $(\partial M\cap (z+ S^{3u}_{r_z}))_{[t]}=\varnothing$ we observe 
	that $((z+ S^{3u}_{r_z})\setminus M)_{[t]} = \varnothing$ and so
	$(t,f_i(t)) \in M^*\cap (z+ S^{3u}_{r_z})$,
	where $i$ is chosen so that $z \in \graph f_i$. 
	If $(\partial M\cap (z+ S^{3u}_{r_z}))_{[t]}\neq \varnothing$   then   
	\eqref{eq:fullProjectionWifhFunctions} follows from
	\eqref{uhvuh} and 
	\eqref{eq:twoInclusionsOne}. 
	Now fix an arbitrary $z=(x,y)\in M^*$
	with $x\in (0,\rho)$ and denote $C\coloneqq  M^* \cap (z+ S^{3u}_{r_z})$. So $C$ is compact and
	thus \eqref{eq:fullProjectionWifhFunctions} implies that we can correctly define 
	$$
	g(t)=\min C_{[t]},\quad t\in (x-r_z,x+r_z).
	$$
	Then $g$ is continuous on $(x-r_z,x+r_z)\cap (0,\rho)$. 
	Indeed,  the compactness of $C$ easily implies that  $g$ is  lower semicontinuous on $(x-r_z,x+r_z)\cap (0,\rho)$.  To prove       moreover, the upper semicontinuity of $g$, consider
	$t\in (x-r_z,x+r_z)\cap (0,\rho)$ and observe that  \eqref{eq:fullProjectionWifhFunctions}
	applied to $z^*\coloneqq  (t,g(t))$ implies that
	\begin{equation*}
	g(\tau)-g(t)\leq 3u|\tau-t|,\quad  \tau \in (t-r_{z^*},t+r_{z^*})\cap (x-r_z,x+r_z),
	\end{equation*}
	which implies that $g$ is upper semicontinuous at $t$. By Lemma \ref{lipmix} we obtain 
	that $g$ is $u$-Lipschitz on $(x-r_z,x+r_z)\cap (0,\rho)$. 
	
	Now, choosing	$\delta>0$ such that
	$[x-\delta,x+\delta] \subset (x-r_z,x+r_z)\cap (0,\rho)$, we obtain the assertion of the lemma.
\end{proof}

Below we will need some easy facts concerning $1$-dimensional DC surfaces in $\R^2$, which are proved
in \cite[Remark 7.1 and Lemma 7.3]{PRZ}. Let us note that in these observations from \cite{PRZ}
the term ``DC graph'' has a different meaning than in \cite{PZ1} and the present article: it 
means there a $1$-dimensional DC surface in $\R^2$.

Thus, in the present terminology, \cite[Remark 7.1]{PRZ} gives the following.
\begin{remark}\label{dcgr2} 
	Let $P\subset \R^2$ be a $1$-dimensional DC surface in $\R^2$
	 and $a \in P$.
	Then
	\begin{enumerate}
		\item[(i)]
		$\Tan(P,a) \cap S^1  $ is a two point set, and
		\item[(ii)]
		there exist $1$-dimensional DC surfaces $P_1, P_2 \subset \R^2$ such that $P \subset P_1 \cup P_2$, $a \in P_1\cap P_2 $
		and $\Tan(P_i,a)$ is a $1$-dimensional space, $i=1,2$.
	\end{enumerate}
\end{remark}

Further, \cite[Lemma 7.3]{PRZ} is the following result.
\begin{lem}\label{zestr}
	Let $P$ be a $1$-dimensional DC surface in $\R^2$ and $0 \in P$.
	Suppose that $\Tan(P,0)$ is a $1$-dimensional space and $(0,1) \notin \Tan(P,0)$. Then there exists 
	$\rho^*>0$ such that, for each $0< \rho< \rho^*$, there exist $\alpha<0 < \beta$ and a DCR function $f$ on $(\alpha, \beta)$ such that
	$P \cap U(0, \rho) = \graph  f|_{(\alpha, \beta)}$.
\end{lem}

 We will also need   the following simple fact which is a standard consequence of the Zorn lemma.

\begin{lem}\label{lipzorn}
	Let $L>0$, $\rho>0$ and $F \subset [0,\rho]\times \R$ be a closed set such that
	\begin{enumerate}
		\item[(i)]
		for each $(x,y) \in F$ with $0<x<\rho$ there exist $\delta>0$ and an $L$-Lipschitz function $g$
		on $[x- \delta, x + \delta]$ such that $g(x)= y$ and $\graph g \subset F$.
	\end{enumerate}
	Then 
	\begin{enumerate}
		\item[(ii)]
		for each $(x,y) \in F$ with $0<x<\rho$ there exists an $L$-Lipschitz function $\gamma$
		on $[0, \rho]$ such that  $\gamma(x)= y$ and 
		$\graph \gamma \subset F$. 
	\end{enumerate}
\end{lem}
\begin{proof}
	To prove (ii), consider an arbitrary $(x,y) \in F$ with $0<x<\rho$. 
	
	Denote by $P$ the set of all $L$-Lipschitz functions $f: (a_f,b_f) \to \R$
	such that  $ (a_f,b_f) \subset (0,\rho)$,    $x \in (a_f,b_f)$ , $f(x)=y$ and  $\graph f \subset F$.
	By (i), we
	obtain $P \neq \varnothing$. Define a partial order on $P$ by inclusion (i.e.,
	$f_1 \leq f_2  \Leftrightarrow \graph f_1 \subset \graph f_2$). Let  $\varnothing \neq T \subset P$
	be a totally ordered set.  Then $\bigcup \{ \graph f: f \in T\}$ is clearly the graph
	of   a function $g \in P$    which is an upper bound of $T$. Consequently, by Zorn theorem,
	$P$ contains a maximal element $f: (a_f,b_f) \to \R$ and we can extend
	$f$  to an $L$-Lipschitz function $\gamma$ on $[a_f,b_f]$.
	Observe that the points
	$(a_f, \gamma(a_f))$, $(b_f, \gamma(b_f))$ belong to $F$
	since the latter set is closed.
	We claim that $a_f=0$ and $b_f=\rho$.
	Indeed, otherwise we can use (i) (applied
	either to $(x,y) = (a_f, \gamma(a_f))$ or  to   $(x,y) = (b_f, \gamma(b_f))$ and easily
	obtain a contradiction with the maximality of $f$.  Consequently, $\gamma$ has all
	 properties from (ii). 
\end{proof}

\begin{lem}\label{zorn}
	Let $M \in \cal D_2$ and $0 \in M$. Let $u>0$ and $A^{4u} \cap \Tan(\partial M,0) \cap S^1 = \{(1,0)\}$. Then there
	exists $r>0$ and an (s)-set $S$ such that
	\begin{equation}\label{parva}
	\partial M \cap A_r^u \subset S \subset M \cap A_r^u. 
	\end{equation}
\end{lem}
\begin{proof}

	By Proposition~\ref{pokr} there exist $\eta>0$ and  $1$-dimensional DC surfaces  
	$P_1,\dots,P_n$  such that 
	$$  \partial M \cap B(0,\eta) \subset  P_1\cup\dots\cup P_n.$$
	Diminishing $\eta$ if necessary, we can suppose that $0\in P_i$ for all $i$.
	Due to Remark~\ref{dcgr2}~(ii) we can  also suppose (changing $n$, if necessary) that $\Tan(P_i,0)$ is a $1$-dimensional linear space for every $i$. 
	Put 
	\begin{equation}\label{eq:definitionOfI}
	 I\coloneqq\{1\leq i\leq n: \Tan(P_i,0)=\spa \{(1,0)\} \}.
	\end{equation} 
	By our assumptions, clearly $I\neq \varnothing$; we can suppose that $I= \{1,\dots,k\}$.
	By our assumptions and the definition of $I$, we can choose $t>0$
	such that  
	$$ \partial M \cap A_t^{3u} \subset  P_1\cup\dots\cup P_k.$$
	Using Lemma \ref{zestr} we obtain that for each $1\leq i \leq k$ there exist 
	$ s_i\in (0,\infty)$ and a DCR function $\vf_i$ on $[0, s_i]$ such that
	$P_i \cap A_{s_i}^{3u} = \graph \vf_i$. 
	Note that (by \eqref{eq:definitionOfI}) $(\vf_i)'_+(0) = 0$ and so, using Lemma \ref{strict}
	we  obtain 
	$0<s\leq \min\{s_1,\dots,s_k\}$  such that, denoting  $f_i \coloneqq  \vf_i |_{[0,s]}$, we have
	that each $f_i$ is $u$-Lipschitz on $[0,s]$ and 
	\begin{equation}\label{jestejednou}
	\partial M\cap A_s^{3u}\subset\bigcup_{i=1}^k  \graph f_i.
	\end{equation}
	Moreover, by the assumptions $0$ is an accumulation point of $\partial M \cap A_s^u$.

	Thus  the assumptions
	of Lemma \ref{L:localExtension} are satisfied. Let $0< \rho< s$ be the corresponding number
	from the assertion of Lemma \ref{L:localExtension}. Set $r\coloneqq  \rho/2$ and
	\begin{equation}\label{defic}
	S\coloneqq   M \cap \bigcup_{i=1}^k \graph f_i \cap ([0,r] \times \R).
	\end{equation}
	Using \eqref{jestejednou}, it is easy to see that \eqref{parva} holds. 
	So it remains to prove that $S$ is an (s)-set. Since $S \subset \bigcup_{i=1}^k \graph f_i $
	and $S \setminus \{0\} \neq \varnothing$, it is sufficient to prove
	that for every  $(x,y) \in S$ with $x\neq 0$ there exists a continuous function
	$h:[0,r] \to \R$  such that $h(x) = y$ and  $\graph h \subset S$. 
	
	To construct $h$, observe that, by the choice of $\rho$, the assertion of 
	Lemma \ref{L:localExtension}  holds. Consequently, for $F\coloneqq    M \cap \bigcup_{i=1}^k \graph f_i \cap ([0,\rho] \times \R)$ and $L\coloneqq  u$, the assumptions of Lemma \ref{lipzorn}
	hold. Therefore there exists a $u$-Lipschitz function $\gamma$ on $[0,\rho]$
	such that $\gamma(x)= y$ and $\graph \gamma \subset F$.

	Consequently the function $h\coloneqq  \gamma|_{[0,r]}$ has  the required properties.
\end{proof}

\begin{lem}\label{L:fillingLemma}
	Let $M$ and $K$ be closed sets in $\er^2$ and let   $x\in M$ and   $\rho>0$ be such that $K\cap U(x,\rho)\subset M$ and $\partial M\cap U(x,\rho)\subset\partial K$.
	If  $d_K$ is DC on $U(x,\rho)$ then $d_M$ is DC on $U(x,\frac{\rho}{2})$.
\end{lem}

\begin{proof}
	 Pick $z\in U(x,\frac{\rho}{2})$.
	First note that if $z\in M$ then $d_M(z)=0$. If $z\notin M$ then  $d_M(z)=d_{\partial M}(z)\geq d_{\partial K}(z)\geq d_K(z)$ since  $d_{M}(z)\leq |x-z|$, $\partial M\cap U(x,\rho)\subset\partial K$ and $\partial K\subset K$. On the other hand $K\cap U(x,\rho)\subset M$ implies  $d_M(z)\leq d_K(z)$ and so $d_M(z)= d_K(z)$.
	Consequently $d_{M}(z)\in\{0,d_K(z)\}$, $z\in U(x,\frac{\rho}{2})     $, and so $d_M$ is DC on $U(x,\frac{\rho}{2})$ by Lemma~\ref{vldc}~(vi).
\end{proof}

Now we can prove the following ``local'' characterisation of $\D_2$ sets by (s)-sets.

\begin{thm}\label{T:localCharacterization}  
	Let $M\subset\er^2$ be a closed set and let $P$ be the set of all isolated points of $M$. Then the following statements are equivalent:
	\begin{enumerate}
		\item\label{cond:dMDC} $M\in \D_2$,
		\item\label{cond:unionOfChSets} for every $z\in \partial M\setminus P$
		there are $\rho>0$, (s)-sets $S_1,\dots,S_m$ and {\it  pairwise different} rotations $\gamma_1,\dots,\gamma_m$ such that
		\begin{equation}\label{lokhr}
		\partial M\cap U(z,\rho) \subset \bigcup_{i=1}^m  (z+ \gamma_i(S_i)) \subset M,
		\end{equation}
		\item\label{cond:unionOfSpecialChSets} for every $z\in \partial M\setminus P$ 
		there are $\rho>0$, (s)-sets $S_1,\dots,S_m$ and  rotations $\gamma_1,\dots,\gamma_m$ such that \eqref{lokhr} holds.
	\end{enumerate}
\end{thm}
%???\begin{rem}
%	If we moreover demand in (2) that the sets  $x+ \gamma_i(S_i)$, $i=1,\dots,m$, are pairwise disjoint,
%	we obtain a condition (2)* which is equivalent to (2).
%\end{rem}
%For nowhere dense $\D_2$ sets teh following following ``local'' and ``global'' characterizations
%follow.

\begin{proof}
	First we prove implication $\eqref{cond:dMDC}\implies\eqref{cond:unionOfChSets}$.
	Let $M$ be a $\D_2$ set and $z\in\partial M\setminus P$. Then $z$ is not an isolated point of $\partial M$ and we obtain by Proposition~\ref{pokr} and Remark~\ref{dcgr2}~(i) that $T\coloneqq\Tan(\partial M,z)\cap S^1$ is a nonempty finite set. Let $T=\{t_1,\dots,t_m\}$ and let $\gamma_i$ be the rotation that maps $(1,0)$ to $t_i$, $i=1,\dots,m$. Since $T$ is finite, there is $u>0$ such that 
	$$
	A^{4u}\cap(\Tan(\gamma^{-1}_i(\partial M-z),0)\cap S^1)=\gamma^{-1}_i(t_i)=(1,0),\quad i=1,\dots,m.
	$$
	By Lemma~\ref{zorn} there is, for every $i=1,\dots,m$, an $r_i>0$ and an (s)-set $S_i$ such that 
	\begin{equation*}
	\gamma^{-1}_i(\partial M-z)\cap A^u_{r_i}=\partial(\gamma^{-1}_i(M-z))\cap A^u_{r_i}\subset S_i\subset \gamma^{-1}_i(M-z),\quad i=1,\dots,m.
	\end{equation*}
	 and consequently
	\begin{equation}\label{inklsi}
	\partial M \cap (z+  \gamma_i(A^u_{r_i})) \subset z + \gamma_i(S_i)
	 \subset M.
	\end{equation}
	By the definition of the tangent cone there is some $\rho>0$ such that
	\begin{equation*}
	\partial M \cap U(z,\rho)\subset\bigcup_{i=1}^m (z+\gamma_i(A^u_{r_i})),
	\end{equation*}
	and so   \eqref{inklsi}  implies that \eqref{lokhr} holds, and the proof of the implication is finished.  
	
Implication $\eqref{cond:unionOfChSets}\implies\eqref{cond:unionOfSpecialChSets}$ is clear and implication $\eqref{cond:unionOfSpecialChSets}\implies\eqref{cond:dMDC}$ follows easily from Corollary~\ref{cor:sSetInD2}, Lemma~\ref{L:fillingLemma} (with $K=\bigcup_{i=1}^m (z+\gamma_i(A^u_\rho))$), \eqref{sjedno} and Lemma~\ref{vldc}~(ii), which concludes the proof of the theorem.
\end{proof}
\begin{rem}\label{simala}
 In (ii) (and in (iii)) we can demand that both $\rho$
 and
$\diam(\bigcup_{i=1}^m (z+\gamma_i(S_i)))$ ``are arbitrarily small'' (i.e., 
smaller than any $\ep>0$    prescribed together with $z\in \partial M \setminus P$). To see this, choose a sufficiently small $0<\rho^*<
 \rho$ and observe that \eqref{lokhr} remains hold if we write
 $\rho^*$ instead of $\rho$ and $S^*_i \coloneqq S_i \cap ([0,\rho^*]
 \times \R)$ (which is an (s)-set by Remark  
\ref{rem:changingRforSset}) instead of $S_i$, $i=1,\dots,m$. 
\end{rem}
\begin{cor}\label{P}
If  $M \in \cal D_2$, then the set $P$ of all isolated points of $M$ is discrete. 
  \end{cor}
  \begin{proof}
 Suppose to the contrary that there exists a point $z\in \overline P \setminus P$. Then clearly $z \in \partial M \setminus P$
  and we can choose $\rho$, $S_1,\dots,S_m$ and $\gamma_1,\dots, \gamma_m$ as in Theorem \ref{T:localCharacterization} (ii);
   so \eqref{lokhr} holds. Since $P \subset \partial M$, there exists $p \in P$ such that
   $p \in  \bigcup_{i=1}^m  (z+ \gamma_i(S_i)) \subset M$, which contradicts the connectivity of $\bigcup_{i=1}^m  (z+ \gamma_i(S_i))$ (cf. Remark \ref{rem:sSetBasicFacts}).  
   \end{proof}

Next we prove Theorem \ref{T:obar} which gives a ``global'' characterization of general $\D_2$ sets by nowhere dense
$\D_2$ sets. 
We will need the following simple observation. 
\begin{lem}\label{L:obar} 
	Let $K$, $M$, $K \subset M$, be closed subsets of $\R^2$. Then the following conditions
	are equivalent
	\begin{enumerate}
		\item\label{cond:systemOfcomponents}
		$M= K \cup C$, where  $C$ is the union of a system of components
		of $\R^2 \setminus K$,
		\item\label{cond:boundaryInclusion}
		$\partial M \subset \partial K.$
	\end{enumerate}
	If these conditions hold and $K \in \D_2$, then  $M \in \D_2$.
\end{lem}

\begin{proof}
  Let \eqref{cond:systemOfcomponents}  hold. Since $C \subset \INt M$, we have
$ \partial M =  M \setminus  \INt M \subset M\setminus C \subset K.$
Obviously, $\partial M \subset \R^2 \setminus \INt K$, and consequently $\partial M \subset K \setminus \INt K = \partial K$.
 We have proved $\eqref{cond:systemOfcomponents}  \Rightarrow \eqref{cond:boundaryInclusion} $.

	To prove \eqref{cond:systemOfcomponents} from \eqref{cond:boundaryInclusion},  it is sufficient to prove that if $D$ is a component of $\er^2\setminus K$ then either $D\subset M $ or $D\cap M=\varnothing$, but it follows from the fact that $\partial M\cap D\subset\partial K\cap D=\varnothing$ and so $M\cap D$ is both open and closed in open connected set $D$.

	To prove the last part of the lemma, it is sufficient to observe that  if $K \in \D_2$  and $M \neq \varnothing$ then  \eqref{cond:boundaryInclusion} together with Lemma~\ref{L:fillingLemma}
	and \eqref{loksem}
	imply that $d_{M}$ is locally DC. Indeed, then $d_M$ is DC on $\er^2$  by Lemma~\ref{vldc}~(ii), and so $M\in\D_2$. 
\end{proof}

\begin{thm}\label{T:obar}
	Let $M\subset\er^2$ be a closed set. Then the following conditions are equivalent.
	\begin{enumerate}
		\item\label{cond:MinD2Obar}  $M \in \D_2$,
		\item\label{cond:existsKNowhereDense}  there exists a nowhere dense $K \in \D_2$  such that $\partial M \subset K\subset M$,
		\item\label{cond:existsKColoring} there exists a nowhere dense $K \in \D_2$  such that  $M= K \cup C$, where  $C$ is the union of a system of components
		of $\R^2 \setminus K$.
	\end{enumerate}
\end{thm}

\begin{proof}
	Implication $\eqref{cond:MinD2Obar}\implies \eqref{cond:existsKNowhereDense}$ follows from Theorem~\ref{T:localCharacterization} as follows.
	
	Suppose that $M\in\D_2$ and let $P$ be the set of all isolated points of $M$.
%	Note that (by Theorem~\ref{T:localCharacterization}~\eqref{cond:unionOfSpecialChSets} and Remark~\ref{rem:sSetBasicFacts}) %$P$ is discrete.
By Corollary \ref{P}, $P$ is discrete. 
	We know that, for every $z\in \partial M \setminus P$, there are $\rho^z>0$, $m_z\in\en$, (s)-sets $S_1^z,\dots,S_{m_{z}}^z$ and  rotations $\gamma_1^z,\dots,\gamma_{m_z}^z$ as in Theorem~\ref{T:localCharacterization}~\eqref{cond:unionOfSpecialChSets}.
	By Remark~\ref{simala} we can suppose $\rho^z\leq 1$ and $\diam S_i^z\leq 1$,  $i=1,\dots,m_z$.  
	The system $\{U(z,\rho^z): z\in \partial M\setminus P \}$ is an open cover of $\partial M\setminus P$. Hence, since $\partial M\setminus P$ is closed and locally compact, we can find $I\subset\en$ such that for every $n\in I$ there are $z_n\in \partial M\setminus P$, $\rho_n>0$, $k_n\in\en$, (s)-sets $S^n_{k}$, $k=1,\dots,k_n$, and isometries $\gamma^n_{k}$, $k=1,\dots,k_n$, such that
	\begin{enumerate}[label=(\alph{*}), ref=\alph{*}]
		\item\label{cond:ballsCover} $\partial M\setminus P\subset\bigcup_{n\in I} U(z_n,\rho_n)$,
		\item\label{cond:ballsLocallyFinite} the system $\{U(z_n,\rho_n):n\in I \}$ is locally finite,
		\item\label{cond:diameterLessThanOne} $\diam S^n_k\leq 1$, $n\in I$, $k=1,\dots,k_n$,
		\item\label{cond:TheInclusions} $\displaystyle
		\partial M\cap U(z_n,\rho_n) \subset \bigcup_{k=1}^{k_n}  (z_n+ \gamma^n_k(S^n_k)) \subset M
		$
		for every $n\in I$.
	\end{enumerate}
	
	Put 
	\begin{equation}\label{defk}
	K\coloneqq P\cup \bigcup_{n\in I}\bigcup_{k=1}^{k_n} (z_n+ \gamma^n_k(S^n_k)).
	\end{equation}
	
	By \eqref{cond:ballsLocallyFinite} and \eqref{cond:diameterLessThanOne} we  obtain  that the system $\{\bigcup_{k=1}^{k_n} (z_n+ \gamma^n_k(S^n_k)):n\in I\}$ is a locally finite system of closed nowhere dense  sets and therefore $K$ is closed nowhere dense.
	Moreover, $\partial M\subset K=\partial K$ by \eqref{cond:ballsCover} and \eqref{cond:TheInclusions} and $K\subset M$ by \eqref{cond:TheInclusions}.
	Finally, $K\in\D_2$ by Corollary~\ref{cor:sSetInD2}, Remark~\ref{rem:satbilityOfD2}~\eqref{cond:D2locallyFiniteStable} and Corollary~\ref{P}.
	
	Implications $\eqref{cond:existsKNowhereDense}\implies \eqref{cond:existsKColoring}$ and $\eqref{cond:existsKColoring}\implies \eqref{cond:MinD2Obar}$ follow from Lemma~\ref{L:obar}.
\end{proof}

\begin{rem}
	\begin{enumerate}
		\item[(i)]
		Lemma \ref{L:obar} shows that each nowhere dense $\cal D_2$ set $K$ yields  (via Lemma \ref{L:obar} (i)) 
		some (sometimes infinitely many) $\cal D_2$ sets $M$ with nonempty interior.
		\item[(ii)]
		A problem whether a given closed set $M \subset \R^2$ belongs to $\cal D_2$ does not reduce
		(by our results) to a problem, whether a corresponding nowhere dense set $K \subset \R^2$ 
		belongs to $\cal D_2$ since there are usually many nowhere dense sets $K\subset M$  with
		$\partial M \subset K$.   Note that these conditions hold for $K\coloneqq  \partial M$,  but \cite[Example 4.1]{PZ1} (or Example	 \ref{nespkomp} below)
		gives an example of $M\in \D_2$ with $\partial M \notin \D_2$.
		
	\end{enumerate}
\end{rem}

Finally, as a consequence of Theorem~\ref{T:localCharacterization} and the proof of Theorem \ref{T:obar}, we easily obtain the following   characterizations  of nowhere dense sets in $\D_2$: 

\begin{thm}\label{charrid}
	Let $M\subset\er^2$ be a nowhere dense closed set and let $P$ be the set of all isolated points of $M$. Then the following conditions are equivalent:
	\begin{enumerate}
		\item\label{cond:MInD2ND} $M\in \D_2$,
		\item\label{cond:localConditionNDspecial}  for every $z\in M\setminus P$ there are $\rho>0$, finitely many (s)-sets $S_1,\dots,S_m$ and {\it  pairwise different} rotations $\gamma_1,\dots,\gamma_m$ such that
		\begin{equation}\label{lokrid}
		M\cap U(z,\rho) = \bigcup_{i=1}^m  (z+ \gamma_i(S_i)) \cap U(z,\rho),
		\end{equation}
		\item\label{cond:localConditionND} for every $z\in M\setminus P$ there are $\rho>0$, finitely many (s)-sets $S_1,\dots,S_m$ and  rotations $\gamma_1,\dots,\gamma_m$ such that \eqref{lokrid} holds,
		\item\label{condd:locallyFiniteUnion}
		there exists a system $(S_{\alpha})_{\alpha\in A}$ of (s)-sets and a system
		$(\gamma_{\alpha})_{\alpha\in A}$ of isometries of $\er^2$ such that the system $(\gamma_{\alpha} (S_{\alpha}))_{\alpha\in A}$ is locally finite and such that
		$$M= P\cup \bigcup_{\alpha\in A}  \gamma_{\alpha} (S_{\alpha}).$$
	\end{enumerate}
\end{thm}
\begin{proof}
Denote by (ii)* (resp. (iii)*) the condition which we obtain if we replace in (ii) (resp. (iii))  equation   \eqref{lokrid}
 by the inclusions
\begin{equation}\label{lokridAlt}
	M\cap U(z,\rho) \subset \bigcup_{i=1}^m  (z+ \gamma_i(S_i)) \subset M.
	\end{equation}
	Since $M = \partial M$,
	 condition  \eqref{lokridAlt} is equivalent to \eqref{lokhr} and so
	the equivalence of (i), (ii)* and (iii)* follows immediately from Theorem~\ref{T:localCharacterization}. Further, \eqref{lokridAlt} clearly implies \eqref{lokrid}, and
	 consequently (ii)* implies (ii) and  (iii)* implies (iii). 
	
	Now we will show $(iii) \Rightarrow (iii)^*$. So suppose that
	(iii) holds and  $x \in M \setminus P$ is given. Find $S_1,\dots, S_m$ and $\gamma_1,\dots, \gamma_m$ by (iii) and choose 
	 $\tilde{\rho}>0$ so small that $\tilde S_i\coloneqq S_i\cap ([0,\tilde{\rho}]\times \er)\subset U(0,\rho)$, $i=1,\dots,m$. Then each $\tilde S_i$ is an (s)-set by Remark~\ref{rem:changingRforSset} and clearly
	$ M\cap U(z,\tilde \rho) \subset \bigcup_{i=1}^m  (z+ \gamma_i(\tilde S_i)) \subset M,$
	 and thus we have proved (iii)*. The above argument proves also
	  $(ii) \Rightarrow (ii)^*$.

Thus we obtain the equivalence of (i), (ii) and (iii).

 To prove $(i) \implies (iv)$, suppose $M \in \cal D_2$ and
  observe that $K$ from Theorem \ref{T:obar} equals to $M$ (by Theorem \ref{T:obar} (ii)). 
   So, chosing $z_n$, $\gamma_k^n$ and $S_k^n$ as in the proof  of Theorem \ref{T:obar}, we obtain that \eqref{defk} holds
    and $M=K$. Since we know that the system of all sets of the form  $z_n+ \gamma_k^n(S_k^n)$ from \eqref{defk}
     is locally finite, (iv) holds.
Finally, implication $\eqref{condd:locallyFiniteUnion}\implies\eqref{cond:MInD2ND}$
follows from Corollary~\ref{cor:sSetInD2}, Corollary \ref{P} and      Remark~\ref{rem:satbilityOfD2}.
\end{proof}

\begin{remark}\label{simala2}
In Theorem \ref{charrid}, $M = \partial M$ and so 
 \eqref{lokhr} implies \eqref{lokrid}. Consequently Remark \ref{simala} shows that, in conditions (ii) and (iii) of Theorem \ref{charrid}, we can demand that   both $\rho$
 and
$\diam(\bigcup_{i=1}^m (z+\gamma_i(S_i)))$ ``are arbitrarily small''. 
\end{remark}

An immediate consequence of Theorem \ref{charrid} is the following 
 result.
\begin{cor}\label{charcomp}
A nonempty nowhere dense perfect compact set is a  $\cal D_2$ set
 if and only if it is a finite union of isometric copies of (s)-sets.
\end{cor}

The following result shows that, in some sense, it suffices to investigate
connected $\cal D_2$ sets only.

\begin{thm}\label{comp}
	A closed set $\varnothing \neq M \subset \R^2$ is a  $\cal D_2$ set
	if and only if
	\begin{enumerate}
		\item[(i)]
		each component of $M$ is a  $\cal D_2$ set
	\end{enumerate}
	and
	\begin{enumerate}
		\item[(ii)]
		the system of all components of $M$ is discrete.
		%system (i.e., for
		%each $z\in \R^2$ there exists $r>0$ such that $U(z,r)$
		%intersects at most one component of $M$).
	\end{enumerate}
\end{thm}
\begin{proof}
Suppose $M \in \cal D_2$ and consider an arbitrary $z\in M$ and the component $C_z$ of $M$ that contains $z$. To prove (ii), we will find $\rho>0$ such that 
\begin{equation}\label{jedncomp}
\text{$C_z$ is the only component of $M$
 intersecting $U(z,\rho)$.}
\end{equation}
 The existence of $\rho$ is obvious if $z$ is an isolated point
 of $M$. Otherwise we can find by Theorem \ref{T:localCharacterization}     
 $\rho>0$, (s)-sets $S_1,\dots,S_m$ and  rotations $\gamma_1,\dots,\gamma_m$ such that \eqref{lokhr} holds. Using Remark 
\ref{rem:sSetBasicFacts}, we obtain 
\begin{equation}\label{sjcomp}
\bigcup_{i=1}^m (z+ \gamma_i(S_i)) \subset C_z.
\end{equation}
Let $C$ be a component of $M$ with $C\cap U(z,\rho) \neq \varnothing$.

If $\partial C\cap U(z,\rho) = \varnothing$, then 
$C\cap U(z,\rho)$ is a nonempty and both open and closed in $U(x,\rho)$, so $U(z,\rho) \subset C$ and thus $C=C_z$.

If  $\partial C\cap U(z,\rho)\neq  \varnothing$, choose a point
$q \in \partial C\cap U(z,\rho)$ and observe that $q \in \partial C \subset C$. Further, since $q \in \partial C \subset \partial M$, 
 by \eqref{lokhr} and \eqref{sjcomp} we obtain   $z\in C_z$. Thus
 $C=C_z$ and (ii) is proved.

To prove (i), consider a component $C$ of $M$. To prove $C\in \cal D_2$, by Lemma \ref{vldc} (ii) it is sufficient to show that $d_C$ is locally DC on $\R^2$. Using \eqref{loksem}, we see that it is sufficient to show that, for each $z\in C$, the function $d_C$ is DC on a
 neigbourhood of $z$. By (ii) we can choose $\rho>0$ such that
 \eqref{jedncomp} holds. Then clearly  $d_C = d_M$ on $U(z,\rho/2)$,
 and thus $d_C$ is DC on  $U(z,\rho/2)$. 

Finally note that if (i) and (ii) hold, then $M\in \cal D_2$ by
    Remark \ref{rem:satbilityOfD2} (ii).   
\end{proof}

\section{Properties of $\cal D_2$ sets and images of $\cal D_2$ sets}

First we will prove several properties of (s)-sets. Then, using our characterization theorems, we will
 obtain some results on general $\D_2$ sets. Finally we will prove
 Theorem  \ref{obrd2} on the stability of $\D_2$ sets with respect to some deformations.

Recall that we already  mentioned  some simple properties
 of (s)-sets; see Remark~\ref{rem:sSetBasicFacts} and  Remark~\ref{rem:changingRforSset}. 

Further mention that \eqref{Ch1} easily implies that, for each (s)-set $S$,
\begin{equation}\label{tans}
\Tan(S,(0,0)) \cap S^1= (1,0).
\end{equation}
An easy consequence of ``mixing lemmas'' is the following fact.

\begin{lem}\label{vlh}
Let $S \subset \R^2$ be an (s)-set and $\pi_1(S)\eqqcolon [0,r]$. Then there exists $K>0$
 such that each continuous $f: [0,r] \to \R$ with $\graph f \subset S$ is a $K$-Lipschitz
 DCR function.
\end{lem}
\begin{proof}
Let $f_1,\dots, f_k$ and $H$ be as in Definition \ref{Ch}.
By Lemma \ref{dcr} (ii) we can choose $K>0$ such that all $f_i$, $1\leq i \leq k$, are $K$-Lipschitz
 functions. Let  $f: [0,r] \to \R$ be a continuous function with $\graph f \subset S$. Then,
 using \eqref{Ch1}, we obtain that
 $f$ is
 $K$-Lipschitz by Lemma \ref{lipmix} and  DCR by Lemma~\ref{L:DCRmixing}.
\end{proof}
\begin{cor}\label{vlash}
If $S$ is an (s)-set and $H$ is as in \eqref{Ch2}, then all functions $h \in H$ are DCR functions,
 and they are equally Lipschitz.
\end{cor}

\begin{lem}\label{hspoc} 
Let $S \subset \R^2$ be an (s)-set with $\pi_1(S)\eqqcolon [0,r]$ and let $H$ be as in 
 \eqref{Ch2}. Then there exists a countable set $H^* \subset H$ such that
 $S= \bigcup_{h\in H^*}\graph h$. 
\end{lem}
\begin{proof}
By \eqref{Ch2}, we have
\begin{equation}\label{ha}
S= \bigcup_{h \in H} \graph h
\end{equation}
 and, by Corollary \ref{vlash}, there exists $K>0$ such that each function $h \in H$ is $K$-Lipschitz.
 Further choose $k \in \en$ by \eqref{Ch1}.

Now choose (using the definition of $k$ and  \eqref{ha}), 
 for each $t \in \Q\cap [0,r]$,  
 functions  $h_1^t,\dots, h_k^t \in H$ such that
$$  S_{[t]} = \{ h_1^t (t),\dots, h_k^t (t)\}$$
 and set
$$ H^*\coloneqq  \bigcup \{ h_i^t:\ t \in \Q\cap [0,r],\ 1 \leq i \leq k\}.$$
Then $H^*$ is countable. To prove 
\begin{equation}\label{hhve}
S= \bigcup_{h \in H^*} \graph h,
\end{equation}
 consider an arbitrary point
 $(x,y) \in S$. Since $S_{[x]}$ is finite, we can choose $\ep >0$ such that
 $ S_{[x]}   \cap (y-\ep, y+ \ep) = \{y\}$. Further choose $x^* \in \Q\cap [0,r]$ such that
 $|x-x^*|< \ep (2K)^{-1}$. By \eqref{ha} there exists $h \in H$
 with $h(x) = y $. Since $h(x^*) \in S_{[x^*]}$,  by the definition of $H^*$ there exists $h^* \in H^*$ with $h^*(x^*) = h(x^*)$. Since $h$, $h^*$ are $K$-Lipschitz, we have
$$ |h(x^*)- h(x)| \leq K|x-x^*| < \frac{\ep}{2},\ \ |h(x^*)- h^*(x)|=
|h^*(x^*)- h^*(x)| \leq K|x-x^*| < \frac{\ep}{2}$$
 and so  $|h^*(x)-y| = |h^*(x)-h(x)|< \ep$. Since $h^*(x) \in S_{[x]}$, we have $h^*(x) =y$
 and thus \eqref{hhve} follows.
\end{proof}
Corollary \ref{vlash} and Lemma \ref{hspoc} have the following immediate consequence.
\begin{cor}\label{sspocgr}
Each (s)-set is a countable union of DC graphs.
\end{cor}
 The following result easily follows.
\begin{prop}\label{sjdcgr}
Each nowhere dense $\cal D_2$ set $M$ is a countable union of DC graphs.
\end{prop}
\begin{proof}
The statement follows from Theorem \ref{charrid} (iv), Corollary \ref{P}, Corollary \ref{sspocgr} and the easy fact that the image
 of a DC graph under an isometry of $\R^2$ is a DC graph.
\end{proof}

\begin{prop}\label{souv}
\begin{enumerate}
 \item[(i)] Each $\cal D_2$ set $M$ is locally pathwise connected; in particular it is 
 locally connected.
\item[(ii)] Each connected $\cal D_2$ set $M$ is pathwise connected.
 Moreover, any two points $x, y \in M$ can be connected by
 a rectifiable curve  lying in $M$.
\end{enumerate}
\end{prop}
\begin{proof}
Let $z\in M$, $r>0$, and $U\coloneqq  U(z,r)\cap M$. To prove (i), it is sufficient to find
 a pathwise connected neigbhourhood  $V\subset U$ of $z$ in the subspace $M$. If $z$
 is an isolated or an interior point of $M$, the existence of $V$ is obvious. Otherwise
 $z \in \partial M$ and we can find by Remark \ref{simala} $\rho \in (0,r)$, (s)-sets
 $S_1,\dots, S_m$ and rotations $\gamma_1,\dots, \gamma_m$ such that
\begin{equation}\label{dveinkl}
\partial M \cap U(z,\rho) \subset Z \subset M\ \ \ \text{and}\ \ \ Z \subset U(z,r),
\end{equation}
  where  $Z\coloneqq  \bigcup_{i=1}^m (z+ \gamma_i(S_i))$.
 Set 
 \begin{equation}\label{defV}
 V\coloneqq  (M\cap U(z,\rho)) \cup Z.
 \end{equation}
  It is clearly sufficient to prove that $V$ is
pathwise connected. Note that Remark \ref{rem:sSetBasicFacts} implies that $Z$ is  pathwise connected   and consider
 an arbitrary $y \in (M\cap U(z,\rho))\setminus Z$. Using 
\eqref{dveinkl}, we obtain
 $y \in \inter M \cap U(z,\rho)$. It is easy to show that there exists a point
 $w \in  \overline{y,z} \cap \partial M$ such that $\overline{y,w} \subset M$. Since
 clearly $\overline{y,w} \subset U(z,\rho)$, we have $\overline{y,w} \subset V$.
 So $y$ can be connected by a path in $V$ with the point $w$ which belongs  to $Z$ by 
\eqref{dveinkl}. Consequently, $V$ is
 pathwise connected.
 
 The first part of (ii) holds since every connected, locally pathwise connected topological space is pathwise connected
  (see, e.g., \cite[Theorem 27.5]{Wi}). To argue that the ``moreover part'' holds, we will say (for a while) that
   a set $A\subset \R^2$ is {\it r-path connected}, if any two points $x$ and $y$ in $A$ can be connected in  $A$
    by a rectifiable path. Corollary \ref{vlash} implies that each (s)-set is r-path connected.  Consequently the argument in the proof of (i) gives that  each $V$ as in \eqref{defV} is even  r-path connected (and thus $M$ is ``locally r-pathwise connected''). So an obvious modification of the
     (standard easy) proof of \cite[Theorem 27.5]{Wi} gives that $M$ is  r-path connected.
     \end{proof}
 \begin{rem}\label{turn}
 Using a straightforward easy (but not trivial) modification of the proof of  \cite[Theorem 27.5]{Wi}, we can obtain the following stronger result:
 
 If $M$ is a connected $\D_2$ set and $x\neq y \in M$, then there exist numbers $t_1< t_2<\dots<t_m$ and a continuous injective
  $f: [t_1,t_m] \to M$ such that  $f(t_1)=x$, $f(t_m) = y$ and each set $f([t_k,t_{k+1}])$, $k=1,\dots, m-1$, is
   a DC graph.
   
   (Note that this statement is equivalent to the assertion that every  $x\neq y \in M$ can be connected in $M$ by a simple curve
    of finite turn; for the notion of the turn see, e.g., \cite{Duda}.)
    
    Indeed, it is not difficult to see that each (s)-set and consequently also each  $V$ as in \eqref{defV} has this
     connectivity property.
 \end{rem}

 For each $\cal D_2$ set $M$, the system of all components of
 $M$ is  discrete (and so countable)  by Theorem \ref{comp}. In the following example we show that the system of all components of 
$\partial M$
 can be uncountable.
\begin{example}\label{nespkomp}
	Let $C\subset [0,1]$ be the classical Cantor ternary set and let 
	 $\{I_n: n \in \en\}$ be all bounded components of $\R \setminus C$.
	For each $n \in \en$, choose an interval $[u_n,v_n] \subset I_n$
	 and set $F\coloneqq  \overline{ \bigcup_{n\in \en} [u_n,v_n]}$.
	 Then $f\coloneqq  (d_F)^2$ is DC on $\R$  (see, e.g., \cite[p. 976]{BB})
	  and so the set
		\begin{equation*}
	K\coloneqq \graph f\cup\graph(-f)
	\end{equation*}
	is a (nowhere dense) $\D_2$ set by \eqref{hlpz1} and \eqref{sjedno}.
	Put
	\begin{equation*}
	M\coloneqq    \{(x,y):\ y\geq f(x)\} \cup  \{(x,y):\ y\leq -f(x)\}.
	\end{equation*}
	 Then $M$ is a $\D_2$ set by Lemma \ref{L:obar}.
	It is easy to see that   $\pi_1(\partial M) = \R \setminus 
	\bigcup_{n\in \en}  (u_n,v_n)$ and so  $\pi_1(\partial M)$
			  has uncountably many components. Consequently
				$\partial M$ has uncountably many components as well.
				
	(In particular, $\partial M$ is not  a $\D_2$ set by 
				Theorem \ref{comp}.)  
\end{example}

We already observed (see  Corollary \ref{P}) that, for each $\cal D_2$ set $M$, the set of all isolated points of
 $M$ is discrete. Now we prove a related result concerning  exceptional points of $\D_2$ sets
 of another type.
 
\begin{prop}\label{tang}
Let $M$ be a  $\cal D_2$ set. Then the set
\begin{equation}\label{defem}
E_M\coloneqq  \{ z \in M:\  \card(\Tan(M,z) \cap S^1) = 1\}
\end{equation}
is discrete.
\end{prop}
\begin{proof}
First consider the case when $M$ is an (s)-set; let $r>0$, $f_1,\dots,f_k$ and $H$ be
 as in Definition \ref{Ch}. Then
\begin{equation}\label{pros}
E_M \cap U(0,r) \subset \{0\}.
\end{equation}
Indeed, if $z=(x,y) \in (M \cap U(0,r)) \setminus \{0\}$, then $0<x<r$ and by 
 \eqref{Ch2} there exists $h \in H$ with $h(x)=y$. Since $h$ is a DCR function by
 Corollary \ref{vlash}, we have $z \notin E_M$ (e.g., by Remark \ref{dcgr2} (i)). 

Further consider the case when $M$ is a nowhere dense $\D_2$ set.
Let $P$ be the set of all isolated points of $M$. Obviously,
  for each $z \in (\R^2\setminus M) \cup P$ there is an $\omega>0$ such that
 $E_M \cap U(z,\omega)= \varnothing$.  If $z \in M\setminus P$, let $\rho>0$, $S_1,\dots,S_m$
 and $\gamma_1,\dots, \gamma_m$ be as in Theorem \ref{charrid} (iii).
 Using \eqref{pros} for $M=S_i,\ i=1,\dots,m,$ we easily obtain $\omega>0$ such that 
 $E_M \cap U(z,\omega)\subset \{z\}$ and conclude that $E_M$ is a discrete set.

Finally consider the case of a general $\D_2$ set $M$. Let $M= K \cup C$ be the decomposition
 of $M$ from   Theorem \ref{T:obar} (iii). Since $K$ is a nowhere dense  $\D_2$ set, we know that $E_K$ (defined as in
 \eqref{defem}) is a discrete set. Thus it is sufficient to prove $E_M \subset E_K$.
 To this end, consider an arbitrary point $z \in E_M$. Then clearly $z \notin C$ and consequently
 $z \in K$. It is easy to see that $z$ is not an isolated point of $K$ and therefore
   $\card(\Tan(K,z) \cap S^1) \geq 1$. Since $\Tan(K,z) \subset \Tan(M,z)$,
	 we obtain $z \in E_K$, which completes the proof.
\end{proof}

\bigskip

An important application of our characterizations of
 $\D_2$ sets is Theorem   \ref{obrd2} below on images of
 $\D_2$ sets. First we prove a lemma on images of (s)-sets.

\begin{lem}\label{obrs}
Let $0\in G \subset \R^2$ be an open set, $c>0$, and let
 $F: G \to \R^2$ be a   locally DC      mapping such that $F(0)=0$ and $F'_+(0, (1,0)) = (c,0)$.
 Let $S \subset G$ be an (s)-set. Then there exist $a>0$ and $b>0$ such that
\begin{equation}\label{shv}
S^*\coloneqq  F\left( S \cap ((-\infty,a] \times \R)\right)  \cap ((-\infty,b] \times \R)
\end{equation}
 is an (s)-set.
\end{lem}
\begin{proof}
  
First note that $F$ is locally Lipschitz on $G$ by Lemma \ref{vldc} (iii).  
Let $f_1,\dots,f_k$ be DCR functions on $[0,r]$ and $H$ be a set of continuous functions
 on $[0,r]$ as in Definition \ref{Ch}. Without any loss of generality, we can suppose
 that $\graph f_i \subset G$, $i=1,\dots,k$. Indeed, otherwise we can diminish $r>0$
 and use Remark \ref{rem:changingRforSset}.

For $i=1,\dots,k$, let
$$ \vf^i(x) = (\vf^i_1(x), \vf^i_2(x))\coloneqq  F((x,f_i(x))),\ \ x\in [0,r].$$
Since $f_i$ is a DCR function, we can find $\ep_i>0$ and a DC extension
 $\tilde f_i: (-\ep_i, r+\ep_i)\to \R$  of $f_i$ such
 that  $(x, \tilde f_i(x)) \in G,\ x \in (-\ep_i, r+\ep_i)$.
Then 
$$\tilde \vf^i(x)=  (\tilde \vf^i_1(x),\tilde  \vf^i_2(x)) \coloneqq  F((x,\tilde f_i(x))),\ \ x\in (-\ep_i, r+\ep_i),$$
 is a DC  mapping   by Lemma \ref{vldc} (ii),(iv). Concequently $\vf^i_1$ and    $\vf^i_2$
 are DCR functions on $[0,r]$ by Lemma \ref{dcr}($(i)\Leftrightarrow (v)$).
Let $\eta_i(x): = (x,\tilde f_i(x)),\ x \in (-\ep_i, r+ \ep_i)$. Then $(\eta_i)'_+(0)= (1,0)$
 and consequently, by the chain rule for one-sided directional derivatives (see, e.g. 
\cite[Proposition 3.6 (i) and Proposition 3.5]{Sh}),
$$ (\tilde \vf^i)'_+(0) = (c,0),\ \ (\vf^i_1)'_+(0)= (\tilde \vf^i_1)'_+(0)=c,\ \    (\vf^i_2)'_+(0)
 = (\tilde \vf^i_2)'_+(0) =     0.$$
Consequently Lemma \ref{strict} gives that
$c$ is
 the strict right derivative of  
 $\tilde \vf^i_1$ at $0$, which easily implies that
there exist  $0<r_i < r$ and $0<\rho_i$   such that $\psi_i\coloneqq \vf^i_1|_{[0,r_i]}$ is an increasing  DCR function and $\psi_i:[0,r_i] \to [0,\rho_i]$ is a bilipschitz bijection.
 Then   Lemma \ref{vldcr} implies that   $h_i\coloneqq  \vf^i_2 \circ (\psi_i)^{-1}$ is a DCR function on $[0,\rho_i]$
  and it is easy to see that
$$  F(\graph (f_i|_{[0,r_i]})) = \graph h_i\ \ \ \text{and}\ \ (h_i)'_+(0) = 0.$$
 Set
$$  a\coloneqq  \min(r_1,\dots,r_k),\  b\coloneqq  \min(\vf_1^1(a),\dots,\vf_1^k(a))\ \  \text{and}\ \ 
f^*_i\coloneqq  h_i|_{[0,b]},\ i=1,\dots,k.$$
Then clearly the set $S^*$ from \eqref{shv} satisfies  $S^* \subset \bigcup_{i=1}^k  \graph f_i^*$.

For each $h\in H$ put
$$  E_h\coloneqq  F(\graph h|_{[0,a]}) \cap ((-\infty,b] \times \R).$$
 Since   $S^* = \bigcup_{h\in H}  E_h$, to prove that $S^*$ is an (s)-set it suffices to show
 that, for each $h\in H$, the set $E_h$ is a graph of a continuous function $h^*$ on $[0,b]$. 
 Since $E_h$ is compact (and each function with compact graph is continuous), it is sufficient to prove that
\begin{equation}\label{jegr}
\text{$E_h$ is a graph of a function $h^*$ on $[0,b]$. }
\end{equation}
Set $\omega(x)\coloneqq  \pi_1(F((x,h(x)))),\ x \in [0,a]$. Then $\omega$ is continuous and consequently
 $\omega([0,a])$ is a closed interval. Since $\omega(0)=0$  and, for some $1\leq i \leq k$,
$ \omega(a)= \vf^i_1(a)\geq b $, we obtain  $[0,b] \subset \omega([0,a])$.
 So, to prove \eqref{jegr}, it is sufficient to show that $\omega$ is injective.
 Suppose, to the contrary, that there exist  $0\leq x_1 < x_2 \leq a$ such that
 $\omega(x_1) = \omega(x_2)$. Set $u_0\coloneqq  x_1$. Further observe that there exists 
 $1 \leq i_0 \leq k$ 
    such that  $h(x_1) = f_{i_0}(x_1)$
  and $u_1\coloneqq  \max\{u \in [x_1,x_2]: h(u) = f_{i_0}(u)\} > u_0$. Then clearly
 either $u_1= x_2$ or we can choose  $1\leq i_1 \leq k$ such that
 $h(u_1) = f_{i_1}(u_1)$ and 
  $u_2\coloneqq  \max\{u \in [x_1,x_2]: h(u) = f_{i_1}(u)\}> u_1$. 
 Proceeding in this way, we obtain numbers  $x_1=u_0 < u_1<...< u_q=x_2$ with $1 \leq q \leq k$ 
 and pairwise different indeces  $i_0, i_1,\dots,i_{q-1}$ such that
 $h(u_k) = f_{i_k}(u_k)$ and $h(u_{k+1}) = f_{i_k}(u_{k+1})$ for each $0\leq k \leq q-1$.
 Then $\omega(u_k) = \vf^{i_k}_1(u_k) < \vf^{i_k}_1(u_{k+1}) = \omega(u_{k+1})$ and therefore
 $\omega(x_1) = \omega(u_0) < \omega(u_1)<\dots< \omega(u_q) = \omega(x_2)$, a contradiction
 which completes the proof.
\end{proof}
\begin{thm}\label{obrd2}
Let $ G \subset \R^2$, $G^*\subset \R^2$ be open sets and let
 $F: G \to G^*$ be a  bijection which is locally bilipschitz and locally DC. Let $M \subset G$ be a $\cal D_2$ set such that $F(M)$ is a closed set.
 Then $F(M)$ is a $\cal D_2$ set.
\end{thm}
\begin{proof}
First consider the case when  $M$ is nowhere dense.

To prove that $M^*\coloneqq  F(M) \in \cal D_2$, we will verify the validity of condition (iii) of Theorem
 \ref{charrid} for $M^*$.  To this end,   consider     an arbitrary point  $z^*\in M^*$ which is not an isolated
 point of $M^*$ and set $z\coloneqq  F^{-1}(z^*)$. Since $M \in \cal D_2$, by Theorem \ref{charrid} (iii)
 there exist $\rho>0$, (s)-sets  $S_1,\dots, S_m$ and rotations $\gamma_1,\dots, \gamma_m$
 such that 
\begin{equation}\label{lokrid2}
M\cap U(z,\rho) = \bigcup_{i=1}^m  (z+ \gamma_i(S_i))\cap U(z,\rho).
\end{equation}

   Remark \ref{simala2} shows that
 we can suppose  that  
\begin{equation}\label{vge}
U(z,\rho) \subset G\ \ \ \text{and}\ \ \ z+ \gamma_i(S_i) \subset G,\
 i=1,\dots,m.
\end{equation}

For each $i=1,\dots,m$, we will apply Lemma \ref{obrs} in the following way. Set $v_i\coloneqq  \gamma_i((1,0))$ and  $w_i\coloneqq  F'_+(z,v_i)$.  Since $F$ is locally bilipschitz, we have $w_i \neq 0$ and consequently we can choose a rotation $\gamma^*_i$ and $c_i>0$ such that 
 $\gamma^*_i((c_i,0)) = w_i$. Now, for each $i=1,\dots,m$, set
$ \beta_i(u)\coloneqq  z+ \gamma_i(u)$, $ \beta_i^*(u) \coloneqq  z^*+ \gamma_i^*(u)$, $u\in \R^2$, and 
\begin{equation}\label{deffi}
F^i\coloneqq  (\beta_i^*)^{-1} \circ F \circ \beta_i.
\end{equation}
 Then $F^i$ is a locally bilipschitz and locally DC bijection from  $G_i\coloneqq  (\beta_i)^{-1}(G)$ onto 
    $G^*_i\coloneqq  (\beta_i^*)^{-1}(G^*)$,    $F^i(0)= 0$ and $(F^i)'(0,(1,0)) = (c_i,0)$. Since $S_i \subset G_i$ by \eqref{vge},   
Lemma \ref{obrs} implies that  there exist $a_i>0$ and $b_i>0$ such that    
$$ S_i^* \coloneqq  F^i(S_i \cap ((-\infty, a_i] \times \R)) \cap
 ((-\infty,b_i] \times \R)$$    
 is an (s)-set. Now choose $\rho^*>0$ so small that
\begin{equation}\label{rhmale}
  \rho^* < \min(b_1,\dots,b_m)\ \ \ \text{and}\ \ \ 
 \diam F^{-1}(U(z^*,\rho^*)) < \min(\rho, a_1,\dots,a_m).
\end{equation}
Set $V\coloneqq  F^{-1}(U(z^*,\rho^*))$.
Then clearly, for each $i$,
$$ (\beta_i)^{-1}(V) \subset (-\infty, a_i) \times \R\ \ \text{and}\ \ (\beta_i^*)^{-1}(U(z^*,\rho^*)) \subset
 (-\infty, b_i) \times \R$$
 and consequently
\begin{equation}\label{rovnajse}
F(\beta_i(S_i) \cap V) = \beta_i^*(S_i^*) \cap U(z^*,\rho^*).
\end{equation}
Since $V\subset U(z,\rho)$ by \eqref{rhmale}, using \eqref{lokrid2} 
 we obtain
 $M\cap V = \bigcup_{i=1}^m (\beta_i(S_i) \cap V)$. Consequently
  \eqref{rovnajse} implies
$$M^* \cap U(z^*,\rho^*) = F(M\cap V) = \bigcup_{i=1}^m 
 \beta_i^*(S_i^*) \cap U(z^*,\rho^*)$$
 and thus condition  Theorem  \ref{charrid} (iii) holds for $M^*$.

To finish the proof, consider an arbitrary $\cal D_2$ set $M\subset G$. By Theorem
 \ref{T:obar} there exists a nowhere dense  $\cal D_2$ set $K$ such that
 $\partial M \subset K \subset M$. Observe that $K^*\coloneqq F(K)$ is closed, since it is clearly closed
 in $G^*$ and $\overline{K^*} \subset M^*\coloneqq F(M) \subset G^*$. Consequently
 $\partial M^* = F(\partial M) \subset F(K)= K^*$ and so  $\partial M^* \subset K^* \subset M^*$.
 Since $K^*$ is clearly nowhere dense and   $K^* \in\cal D_2$  by the first part of the proof,
  Theorem
 \ref{T:obar} implies that $M^* \in\cal D_2$.
\end{proof}

%\section{Properties and examples of (s)-sets???}

\end{document}